\documentclass[12pt]{amsart}

\usepackage{amsmath, amsthm, amssymb, fullpage, amscd}

\usepackage[english]{babel}
\usepackage{eucal}
\usepackage{lineno}

\usepackage{thmtools}
\usepackage{float} 
\usepackage{thm-restate}
\usepackage{tikz-cd, pgf, tikz}
\usepackage{geometry}
\usepackage[colorlinks, bookmarks=true, backref=page]{hyperref}
\hypersetup{
	colorlinks,
	citecolor=blue,
	linkcolor=blue,
	urlcolor=blue}
\usepackage{xcolor}
\usepackage{mathtools, mathrsfs}
\usepackage{graphicx}
\usepackage{caption}
\usepackage{subcaption}
\usepackage{bm, bbm}
\usepackage{comment}
\usepackage{array}
\usepackage{latexsym}
\usepackage{stmaryrd}
\usepackage{commath,enumitem,chngcntr}
\usepackage[titletoc,toc,title]{appendix}
\usepackage{pst-node}
\mathtoolsset{showonlyrefs}
\usepackage{overpic}

\usetikzlibrary{arrows,calc,positioning}

\usepackage{natbib}
\usepackage{overpic}

\newcommand{\dmo}{\DeclareMathOperator}
\dmo{\Id}{Id}
\dmo{\pt}{pt}
\dmo{\ab}{ab}
\dmo{\GL}{GL}
\dmo{\SO}{SO}
\dmo{\SU}{SU}
\dmo{\SL}{SL}
\dmo{\PSL}{PSL}
\dmo{\im}{Im}
\dmo{\re}{Re}
\dmo{\Ad}{Ad}
\dmo{\ad}{ad}
\dmo{\tr}{tr}
\dmo{\Log}{Log}
\dmo{\Homeo}{Homeo}
\dmo{\Vol}{Vol}
\dmo{\supp}{supp}
\dmo{\Isom}{Isom}
\dmo{\Hom}{Hom}
\dmo{\Ext}{Ext}
\dmo{\Tor}{Tor}
\dmo{\PP}{\mathbf{P}}
\dmo{\Area}{Area}
\dmo{\Length}{Length}
\dmo{\inte}{int}
\dmo{\sech}{sech}
\dmo{\diag}{diag}
\dmo{\arcosh}{arcosh}
\dmo{\artanh}{artanh}
\dmo{\Foot}{Foot}
\dmo{\Mod}{Mod}

\newcommand{\C}{\mathbb C}
\newcommand{\R}{\mathbb R}

\newcommand{\Z}{\mathbb Z}
\newcommand{\N}{\mathbb N}
\newcommand{\T}{\mathcal T}

\newcommand{\F}{\mathcal F}

\renewcommand{\H}{\mathbb H}

\renewcommand{\u}{\mathbf{u}}
\newcommand{\n}{\mathbf{n}}
\newcommand{\ww}{\mathbf{w}}
\newcommand{\frs}{\mathfrak{s}}

\newcommand{\hl}{\mathbf{hl}}
\newcommand{\cl}{\mathbf{l}_\C}

\newcommand{\ep}{\epsilon}
\newcommand{\de}{\delta}
\newcommand{\al}{\alpha}

\renewcommand{\tr}{\mathrm {tr}}

\newcommand{\from}{\colon}

\newtheorem{theorem}{Theorem}[section]

\newtheorem{corollary}[theorem]{Corollary}
\newtheorem{proposition}[theorem]{Proposition}
\newtheorem{lemma}[theorem]{Lemma}

\newtheorem*{goodrosegeometry}{Theorem~\ref{thm:good-rose-geometry}}
\newtheorem*{elevation}{Corollary~\ref{elevation}}

\theoremstyle{definition}
\newtheorem{definition}[theorem]{Definition}

\theoremstyle{remark}
\newtheorem{remark}[theorem]{Remark}

\usepackage{lipsum}
\let\OLDthebibliography\thebibliography
\renewcommand\thebibliography[1]{
	\OLDthebibliography{#1}
	\setlength{\parskip}{2pt}
	\setlength{\itemsep}{0pt plus 0.3ex}
}

\numberwithin{equation}{section}

\title{Surface subgroups with prescribed finite quotients}
\author{Jean-Fran\c cois Lafont, Zhenghao Rao, and Hongbin Sun}
\date{}

\address{Department of Mathematics, The Ohio State University, Columbus, OH 43210, USA}
\email{jlafont@math.ohio-state.edu}

\address{Department of Mathematics, Rutgers University--New Brunswick, Piscataway, NJ 08854, USA}
\email{zhenghao.rao@rutgers.edu}

\address{Department of Mathematics, Rutgers University--New Brunswick, Piscataway, NJ 08854, USA}
\email{hongbin.sun@rutgers.edu}

\subjclass[2020]{Primary 57K32; Secondary 20H10, 20E18, 57R22}

\keywords{Surface subgroups, hyperbolic $3$-manifolds,
  Kleinian groups, finite quotients, profinite topology, atoroidal surface bundles}

\begin{document}

\begin{abstract}
    We prove that every finite quotient of a finite covolume Kleinian group $G$ can be realized by a surface subgroup. 
    More precisely, given an epimorphism from $G$ onto a finite group $F$, there exists a surface subgroup that maps surjectively onto $F$.
    As a consequence, the union of the surface subgroups is dense in $G$ with respect to the profinite topology.
    We also apply this result to obtain more examples of atoroidal surface bundles with signature zero in the construction of Kent--Leininger.
\end{abstract}

\maketitle

\vspace{-6mm}
\tableofcontents
\vspace{-6mm}
\section{Introduction}

In the past two decades, one of the major developments on $3$-manifold groups is the virtual compact specialness theorem theorem of Agol~\cite{Agol} for the fundamental groups of closed hyperbolic $3$-manifolds. This result implies Waldhausen's virtual Haken conjecture and Thurston's virtual fibering conjecture for closed hyperbolic $3$-manifolds, and has many further applications to the topology of $3$-manifolds and $3$-manifold groups.

An important consequence of the virtual compact specialness theorem is that hyperbolic $3$-manifold groups are LERF, a strong separability property.
We will not review the precise definition of LERF here, but the LERFness implies that
for any finitely generated proper subgroup $H$ of a hyperbolic $3$-manifold group $G$,
there exists an epimorphism $\rho:G\to F$ to a finite group $F$ such that $\rho(H)\ne F$. In particular, this holds for any surface subgroup $H<G$. Here, a \emph{surface subgroup} is a subgroup isomorphic to the fundamental group of a closed,
connected, orientable surface of genus at least two.

In this paper, we prove a result in a complementary direction. Given a finite quotient $\rho:G \to F$ of a finite-covolume Kleinian group $G$, we find a surface subgroup $H<G$ that surjects onto $F$ via $\rho$. This result implies that the union of surface subgroups is dense in any finite-covolume hyperbolic $3$-orbifold group, with respect to the profinite topology.

\begin{theorem}\label{main}
    Suppose $G$ is a finite covolume Kleinian group, $F$ is a finite group, and $\rho\from G\to F$ is an epimorphism. Then for any subgroup $A<F$, there exists a quasi-Fuchsian surface subgroup $H<G$ such that $\rho(H)=A$.
\end{theorem}

This result is analogous to that of \cite{LM15}, where quasi-Fuchsian surfaces are constructed to represent integral multiples of prescribed rational second homology classes in closed hyperbolic $3$-manifolds.
As an application, we study the homology of immersed surfaces in finite covers. A null-homologous immersed closed oriented
surface need not have null-homologous connected elevations.
However, if the elevation is unique, it is also null-homologous.

\begin{corollary}\label{elevation}
Let $M$ be a finite volume hyperbolic $3$-manifold, and $p:\bar M\rightarrow M$ a finite cover. Then there exists a $\pi_1$-injective immersed closed, connected, 
oriented surface $j: N \rightarrow M$, such that there is a unique elevation $i:\bar N \rightarrow \bar M$. As a result, for this surface $N$, if 
$j_*[N]=0\in H_2(M)$ then $i_*[\bar N]=0\in H_2(\bar M)$.
\end{corollary}

We then apply this result to the atoroidal surface bundles over surfaces
constructed by Kent--Leininger~\cite{KL}.
For every acceptable branched cover in their construction, we obtain infinitely many surface subgroups for which the associated surface bundle has signature zero. This complements the results of Lafont--Miller--Ruffoni~\cite{LMR}; see Section~\ref{sec:KL}.

To prove Theorem \ref{main}, a key step is to prove the following special but essential case: $G$ is torsion-free and $A=F$.

\begin{theorem}\label{main2}
Suppose $M=\H^3/G$ is a hyperbolic $3$-manifold of finite volume for a Kleinian group $G$, $F$ is a finite group, and $\rho\from G\to F$ is an epimorphism. Then there exists a quasi-Fuchsian surface subgroup $H<G$ such that $\rho(H)=F$.
\end{theorem}

Theorem~\ref{main2} holds for both closed and cusped hyperbolic $3$-manifolds. We focus primarily on cusped hyperbolic $3$-manifolds, since the good pants construction is more complicated for them. The same proof readily applies to closed hyperbolic $3$-manifolds, and is technically easier in this case. Hence, we assume that $M$ is a cusped hyperbolic $3$-manifold in most situations.

The bulk of this paper is devoted to proving Theorem \ref{main2}. The proof of Theorem \ref{main} follows from Theorem \ref{main2} and a subgroup separability argument. 
The main ingredient of the proof of Theorem \ref{main2} is the good pants construction initiated by Kahn--Markovi\'c~\cite{KM12b}. The good pants construction is a program that constructs nearly geodesic $\pi_1$-injective closed immersed surfaces in certain Riemannian manifolds by gluing good pants and other good objects. For closed hyperbolic $3$-manifolds, nearly geodesic $\pi_1$-injective closed immersed surfaces were constructed by Kahn--Markovi\'c~\cite{KM12b}. 
Then this construction was generalized to higher-dimensional closed hyperbolic manifolds by Hamenst\"adt~\cite{Ham15} in odd dimensions and Kahn--Rao~\cite{KR26} in even dimensions, and to cusped hyperbolic $3$-manifolds by Kahn--Wright~\cite{KW21}. 
Here a hyperbolic $3$-manifold is called \emph{cusped} if it is complete, non-compact, and has finite volume.

In this paper, we introduce a new nearly geodesic object in hyperbolic $3$-manifolds called a \emph{good rose}. It is a genus-$g$ surface with a single boundary component, obtained by thickening an immersed rose with $2g$ petals in a nearly geodesic way. We first construct a good rose in $M$ so that the image of its fundamental group surjects onto the prescribed finite quotient via $\rho$; this is easier to arrange since the group of a good rose is a free subgroup of $G$. Using the pants cobordism group (see \cite{LM15, Sun22}), we then close up the good rose by attaching good pants and good hamster wheels so that the resulting closed immersed surface is nearly geodesic and hence $\pi_1$-injective.
A related result is due to Masters~\cite{Masters}, who proved that every finitely generated free subgroup of a cocompact right-angled Kleinian reflection group is contained in a surface subgroup.

\vspace{6pt}

\noindent\textbf{Organization.}
In Section~\ref{pre-section}, we review some preliminary results on the good pants construction. In Section~\ref{rose-section}, we introduce the perfect rose as a hyperbolic surface, the good rose as a nearly geodesic immersed surface in a hyperbolic $3$-manifold, and state a key estimate theorem for good roses, whose proof is deferred to Appendix~\ref{appen:good-rose-geometry}. In Section~\ref{sec:construction}, we first construct a good rose in the hyperbolic $3$-manifold $M$ so that $\rho$ is already surjective when restricted to the group of the rose, then close up the good rose to a nearly geodesic closed immersed surface $\Sigma_0\to M$. Theorem~\ref{main2} is also proved in Section~\ref{sec:construction}, modulo the $\pi_1$-injectivity of $\Sigma_0\to M$. We prove that Theorem~\ref{main2} implies Theorem~\ref{main} in Section~\ref{provemain}. In Section~\ref{sec:KL}, we present an application of the main theorem to atoroidal surface bundles constructed by Kent--Leininger. Finally, in Appendix~\ref{incompressible-section}, we prove that $\Sigma_0\to M$ is a $\pi_1$-injective map.

\vspace{6pt}

\noindent\textbf{Computational resource disclosure.}
All mathematical results in this paper were developed by the authors.
The authors used ChatGPT to search the literature, check proofs, and proofread the manuscript. The authors take full responsibility for the content of the paper.

\vspace{6pt}

\noindent\textbf{Acknowledgments.}
The work of Z. R. is partially supported by NSF grant DMS-2603932. The work of H. S. is partially supported by Simons Collaboration Grant 615229.

\section{Preliminaries}\label{pre-section}

In this section, we review some preliminary concepts and results in hyperbolic geometry.

\subsection{Inefficiency theory}

We first recall the inefficiency theory for hyperbolic $3$-manifolds from \cite{LM15}. See also \cite{Sun22,Rao25}, and \cite{KM15} for the $2$-dimensional case. In this section, $M$ is an oriented hyperbolic $3$-manifold.

\subsubsection{Terminology}\label{sub3.1}

We start with the concepts of framed geodesic segments.

\begin{definition}

An \textit{oriented\ $\partial$-framed\ segment} in $M$ is a triple:
\begin{equation}
    \mathfrak{s}=(s,\vec{n}_{ini},\vec{n}_{ter}),
\end{equation}
where $s$ is an oriented immersed compact geodesic segment, $\vec{n}_{ini}$ and $\vec{n}_{ter}$ are two unit normal vectors at the initial endpoint and terminal endpoint of $s$, respectively.

\begin{enumerate}
    \item The \textit{carrier\ segment} is the oriented segment $s$.

\item The \textit{initial\ endpoint} $p_{ini}(\mathfrak{s})$ and the \textit{terminal\ endpoint} $p_{ter}(\mathfrak{s})$ are the initial endpoint and the terminal endpoint of $s$, respectively.

\item The \textit{initial\ framing} $\vec{n}_{ini}(\mathfrak{s})$ and the \textit{terminal\ framing} $\vec{n}_{ter}(\mathfrak{s})$ are the unit normal vectors $\vec{n}_{ini}$ and $\vec{n}_{ter}$.

\item The \textit{initial\ direction} $\vec{t}_{ini}(\mathfrak{s})$ and the \textit{terminal\ direction} $\vec{t}_{ter}(\mathfrak{s})$ are the unit tangent vectors in the direction of $s$ at the initial point and the terminal point, respectively.

\item The \textit{orientation\ reversal} of $\mathfrak{s}$ is defined to be
\begin{equation}
    \mathfrak{\overline{s}}=(\overline{s},\vec{n}_{ter},\vec{n}_{ini}),
\end{equation}
where $\overline{s}$ is the orientation reversal of $s$. The \textit{framing\ flipping} of $\mathfrak{s}$ is defined as
\begin{equation}
    \mathfrak{s}^*=(s,-\vec{n}_{ini},-\vec{n}_{ter}).
\end{equation}
\end{enumerate}

\end{definition}

It follows from the definition that:
\begin{equation}
    \overline{\mathfrak{s}^*}=\overline{\frs}^*.
\end{equation}

\begin{definition}
Suppose that $\mathfrak{s}$ is an oriented $\partial$-framed segment in $M$. We define the \emph{length} of $\mathfrak{s}$, denoted by $l(\mathfrak{s})\in[0,+\infty)$, to be the length of the unframed segment $s$ carrying $\mathfrak{s}$. The \emph{phase} of $\mathfrak{s}$, denoted by $\phi({\mathfrak{s}})\in \R/2\pi\Z$, is defined as the angle from the initial framing $\vec{n}_{ini}$ to the parallel transport of $\vec{n}_{ter}$ to the initial endpoint of $\mathfrak{s}$ via $s$, signed with respect to the normal orientation induced from $\vec{t}_{ini}$ and the orientation of $M$. We define the \emph{complex\ length} of $\mathfrak{s}$ by
\begin{equation}
    \cl(\mathfrak{s})=l(\mathfrak{s})+i\phi(\mathfrak{s})\in \mathbb{C}/2\pi i\mathbb{Z}.
\end{equation}
We observe that the length and the phase remain unchanged under orientation reversal and framing flipping.

For an oriented closed geodesic curve $c$, we can also talk about its \textit{length,\ phase}, or \textit{complex\ length}, by taking an arbitrary unit normal vector $\vec{n}$ at a point $p\in c$, and regarding $c$ as a $\partial$-framed segment obtained by cutting $c$ at $p$ and equipped with framing $\vec{n}$ at both endpoints. If the complex length of $c$ satisfies the inequality
\begin{equation}
    |\cl(c)-2R|<2\ep,
\end{equation}
for some $R,\ep>0$, then we say $c$ is an $(R,\ep)$-$good$ $curve$. Here we mean that there is a lift of $\cl(c)$ in $\C$ to make the inequality hold. In this paper, we always consider a small value for $\ep$, making the lift unique when it exists.
\end{definition}

\begin{definition}
Let $0\leq\de<\frac{\pi}{3}$, $L>0$, and $0<\theta<\pi$ be constants.
\begin{enumerate}
\item An oriented $\partial$-framed segment $\mathfrak{s}$ is said to be $\delta$-\emph{consecutive} to another oriented $\partial$-framed segment $\mathfrak{s}'$, if the following hold:
\begin{itemize}
    \item the terminal endpoint of $\mathfrak{s}$ is the initial endpoint of $\mathfrak{s}'$;
    \item the terminal framing of $\mathfrak{s}$ is $\de$-close to the initial framing of $\mathfrak{s}'$.
\end{itemize}
We say $\mathfrak{s}$ is \textit{consecutive} to $\mathfrak{s}'$ if $\de$ is 0. When $\mathfrak{s}$ is $\delta$-consecutive to $\mathfrak{s}'$, the \textit{bending\ angle} between $\mathfrak{s}$ and $\mathfrak{s}'$ is the angle between the terminal direction of $\mathfrak{s}$ and the initial direction of $\mathfrak{s}'$. This angle falls within the range $[0,\pi]$ and is denoted by $\angle(\mathfrak{s},\mathfrak{s}')$.

\item A $\de$-\textit{consecutive\ chain} of oriented $\partial$-framed segments is a finite sequence $\mathfrak{s}_1\mathfrak{s}_2\cdots\mathfrak{s}_m$ such that each $\mathfrak{s}_i$ is $\de$-consecutive to $\mathfrak{s}_{i+1}$. It is a $\de$-\textit{consecutive\ cycle} if furthermore $\mathfrak{s}_m$ is $\de$-consecutive to $\mathfrak{s}_1$. A $\de$-consecutive chain or cycle $\mathfrak{s}_1\cdots\mathfrak{s}_m$ is called $(L,\theta)$-tame, if the length of each $\mathfrak{s}_i$ is greater than $2L$ and the bending angle between $\mathfrak{s}_i$ and $\mathfrak{s}_{i+1}$ is less than $\theta$.

\item For a $\de$-consecutive chain $\mathfrak{s}_1\cdots\mathfrak{s}_m$, we define its \textit{reduced\ concatenation}, denoted by $\langle\mathfrak{s}_1\mathfrak{s}_2\cdots\mathfrak{s}_m\rangle$,
to be the oriented $\partial$-framed segment as follows. The carrier segment of $\langle\mathfrak{s}_1\mathfrak{s}_2\cdots\mathfrak{s}_m\rangle$ is the geodesic arc which is homotopic to the concatenation of the carrier segments of the $\mathfrak{s}_i$, relative to the initial endpoint of $\mathfrak{s}_1$ and the terminal endpoint of $\mathfrak{s}_m$; the initial framing of $\langle\mathfrak{s}_1\mathfrak{s}_2\cdots\mathfrak{s}_m\rangle$ is the unit normal vector closest to the initial framing of $\mathfrak{s}_1$; the terminal framing of $\langle\mathfrak{s}_1\mathfrak{s}_2\cdots\mathfrak{s}_m\rangle$ is the unit normal vector closest to the terminal framing of $\mathfrak{s}_m$.

\item For a $\de$-consecutive cycle $\mathfrak{s}_1\cdots\mathfrak{s}_m$, whose conjugacy class is loxodromic, we define its \textit{reduced\ cyclic\ concatenation}, denoted by
\begin{equation}
    [\mathfrak{s}_1\mathfrak{s}_2\cdots\mathfrak{s}_m],
\end{equation}
to be the unframed oriented closed geodesic curve freely homotopic to the concatenation of the carrier segments of each $\mathfrak{s}_i$, provided that the corresponding conjugacy class is loxodromic. 
\end{enumerate}
\end{definition}

\subsubsection{Inefficiency of chains and cycles of framed segments} 

For any bending angle $\theta$, we define the \textit{inefficiency} of $\theta$ as 
\begin{equation}\label{eq:angle-ineff}
    I(\theta)=2\ln(\sec(\theta/2)).
\end{equation}
For sufficiently tame approximately consecutive $\partial$-framed segments and cycles, Liu and Markovi\'c estimated their inefficiency in terms of length and phase. We restate the results as follows.

\begin{lemma}[Lemma 4.8 in \cite{LM15}]\label{lem:ineff}
Given any $\delta\geq0$, $\pi>\theta>0$, and $\ep>0$, there exists $L>0$ such that the following hold for every integer $m\geq 2$. 
\begin{enumerate}
    \item If $\mathfrak{s}_1\cdots\mathfrak{s}_m$ is an $(L,\theta)$-tame $\de$-consecutive chain of oriented $\partial$-framed segments, with $\theta_i\in[0,\theta)$ the bending angle between $\mathfrak{s}_i$ and $\mathfrak{s}_{i+1}$, then
\begin{equation}
    |l(\langle\mathfrak{s}_1\cdots\mathfrak{s}_m\rangle)-\sum\limits_{i=1}^{m}l(\mathfrak{s}_i)+\sum\limits_{i=1}^{m-1}I(\theta_i)|<(m-1)\ep
\end{equation}
and
\begin{equation}
    |\phi(\langle\mathfrak{s}_1\cdots\mathfrak{s}_m\rangle)-\sum\limits_{i=1}^{m}\phi(\mathfrak{s}_i)|<(m-1)(\de+\ep),
\end{equation}
where $|\cdot|$ on $\R/2\pi\Z$ is understood as the distance from 0 valued in $[0,\pi]$.

\item If $\mathfrak{s}_1\cdots\mathfrak{s}_m$ is an $(L,\theta)$-tame $\de$-consecutive cycle of oriented $\partial$-framed segments, with $\theta_i\in[0,\theta)$ the bending angle between $\mathfrak{s}_i$ and $\mathfrak{s}_{i+1}$, and $\mathfrak{s}_{m+1}$ equal to $\mathfrak{s}_1$, then
\begin{equation}
    |l([\mathfrak{s}_1\cdots\mathfrak{s}_m])-\sum\limits_{i=1}^{m}l(\mathfrak{s}_i)+\sum\limits_{i=1}^{m}I(\theta_i)|<m\ep
\end{equation}
and
\begin{equation}
    |\phi([\mathfrak{s}_1\cdots\mathfrak{s}_m])-\sum\limits_{i=1}^{m}\phi(\mathfrak{s}_i)|<m(\de+\ep),
\end{equation}
where $|\cdot|$ on $\R/2\pi\Z$ is understood as the distance from 0 valued in $[0,\pi]$.

\end{enumerate}
\end{lemma}

\subsection{Ingredients from Kahn--Wright's construction}

In this subsection,
we briefly recall the construction of surface subgroups by Kahn--Wright in a cusped hyperbolic 3-manifold $M$. They used two model components, perfect pants and perfect hamster wheels, as building blocks, found a large family of good pants and good hamster wheels in $M$, and glued them in a nearly geodesic way to obtain a closed immersed surface in $M$.

For the purpose of this paper, we recall the following concepts from \cite{KM12b, KW21, Rao25}.

\subsubsection{Orthogeodesics and feet}
Suppose that $\al_i\from[0,1]\to M$ are two geodesics in $M$ with distinct images, and $\eta\from [0,1]\to M$ is a geodesic segment. If for each $i=0,1$ there exists $x_i\in[0,1]$ such that $\eta(i)=\al_i(x_i)$ and $\langle\eta'(i),\al'_i(x_i) \rangle=0$, i.e., $\eta$ is orthogonal to $\al_i$, we call $\eta$ an \emph{orthogeodesic} between $\al_0$ and $\al_1$.

If $\al_0$ and $\al_1$ intersect, it may happen that $\eta$ is a constant function. When $\eta$ is nonconstant, we can reparametrize $\eta$ to $\eta\from[0,a]\to M$ with unit speed, and we then call $\eta'(0)$ and $-\eta'(a)$ the \emph{feet} of $\eta$.

\subsubsection{Slow and constant turning normal fields}
Suppose $\gamma$ is a closed geodesic in $M$. Its unit normal bundle $N^1(\gamma)$ is a torsor for $V_\gamma:=\C/(\cl(\gamma)\Z+2\pi i\Z)$; $N^1(\gamma)$ fibers over $\gamma$, and a unit normal field $v$ for $\gamma$ is a section of $N^1(\gamma)$. For any smooth normal field $v$, we have a curve in $V_\gamma$ up to translation, which then yields a \emph{slope} at each point of $\gamma$.
A \emph{constant turning normal field} for $\gamma$ is a smooth unit normal field $v$ with constant slope. This constant slope is given by $(\theta+2k\pi)/B$, where $\cl(\gamma)=B+i\theta$ and $k\in\Z$.

In the case where $|\theta+2k\pi|<\pi$, we call $v$ a \emph{slow and constant turning normal field}.
When $\gamma$ is an $(R,\ep)$-good curve with small $\ep$, such normal fields exist. Moreover, the space of slow and constant turning normal fields for $\gamma$ is a circle, and every such field is determined by its value at an arbitrary point of $\gamma$.

\subsubsection{Pants and half-lengths}\label{subsec:pants}

Let $P_0$ be an oriented topological pants as a manifold with boundary. A \emph{pair of pants} in $M$ is an injective homomorphism $\pi_1(P_0)\to G$, up to conjugacy, such that it sends boundary components of $P_0$ to loxodromic elements in $G$. A pair of pants in $M$ is determined by a continuous map $f\colon P_0\to M$ up to homotopy, and we also call $f$ or $f(P_0)$ a pair of pants. Let $C_i$ with $i=1,2,3$ be the boundary components of $P_0$, and we orient $C_i$ such that $P_0$ is on the left of $C_i$, for $i\in\Z/3\Z$. Then there is a unique oriented closed geodesic $\gamma_i$ in $M$ that is freely homotopic to $f(C_i)$. We can homotope $f$ such that $f$ maps $C_i$ to $\gamma_i$. In the remainder of the paper, when talking about pants in $M$, we will always assume that their boundaries are sent to closed geodesics.

Moreover, given any simple non-separating arc $\al_i$ connecting $C_{i-1}$ and $C_{i+1}$, we may homotope $f$ such that $f$ maps $\al_i$ to an orthogeodesic $\eta_i$ between $\gamma_{i-1}$ and $\gamma_{i+1}$. 
We refer to the $\gamma_i$'s and the $\eta_i$'s as \emph{cuffs} and \emph{short\ orthogeodesics} of the pants $f$, respectively, and the feet of $\eta_{i-1}$ and $\eta_{i+1}$ on $\gamma_i$ are called the \emph{feet} of $f$ on $\gamma_i$. When $\gamma$ is viewed as a cuff, the information of the associated pants is also carried by $\gamma$. 

Cutting along the three short orthogeodesics gives two isometric oriented right-angled hexagons sharing three sides $\eta_i$, where the orientations are determined by the orientations of the $\gamma_i$. Let $\mathbf{u}_i^j$ be the foot of $\eta_j$ on $\gamma_i$ for $i\in\Z/3\Z$ and $j\neq i$. In this case, the complex distance $\mathbf{d}_{\gamma_i}(\mathbf{u}_i^j,\mathbf{u}_i^k)$ from $\mathbf{u}_i^j$ to $\mathbf{u}_i^k$ always follows the orientation of $\gamma_i$ for distinct $i,j,k$. Since the two right-angled hexagons are isometric, we have
\begin{equation*}
    \mathbf{d}_{\gamma_i}(\mathbf{u}_i^{i-1},\mathbf{u}_i^{i+1})
    =\mathbf{d}_{\gamma_i}(\mathbf{u}_i^{i+1},\mathbf{u}_i^{i-1}).
\end{equation*}
Thus we let
\begin{equation*}
    \hl(\gamma_i)=\mathbf{d}_{\gamma_i}(\mathbf{u}_i^{i-1},\mathbf{u}_i^{i+1})
\end{equation*}
be the \emph{half-length} of $\gamma_i$. It is not hard to see $2\hl(\gamma_i)=\cl(\gamma_i)$.
For $R,\ep>0$, a pants $P$ in $M$ is said to be $(R,\ep)$-\emph{good}, if for each cuff $\gamma$, we have
\begin{equation}
    |\hl(\gamma)-R|<\ep.
\end{equation}

\subsubsection{Quotient normal bundles and shears}

Suppose that $\gamma$ is a cuff of some pants in $M$. 
Since $N^1(\gamma)$ is a torsor for $\C/(\cl(\gamma)\Z+2\pi i\Z)$,
there is an involution $\sigma$ on $N^1(\gamma)$ given by $\sigma(x)=x+\hl(\gamma)$.
We denote by $\N^1(\sqrt{\gamma})$ the quotient of $N^1(\gamma)$ by $\sigma$, and call it the \emph{quotient unit normal bundle} of $\gamma$, which is a torsor for $\C/(\hl(\gamma)\Z+2\pi i\Z)$. When $\gamma$ is a cuff of a pair of pants $\Pi$, by the discussion in Section \ref{subsec:pants}, the two feet of $\Pi$ on $\gamma$ descend to the same point in $N^1(\sqrt{\gamma})$. We denote by $\mathbf{foot}_\gamma(\Pi)$ their image in $N^1(\sqrt{\gamma})$ and call it the foot of $\Pi$ on $\gamma$ in $N^1(\sqrt{\gamma})$.

For an oriented closed geodesic $\gamma$ and two unit normal vectors $\u_0,\u_1$, we define the shear from $\u_0$ to $\u_1$ along $\gamma$ to be $$\u_1-\u_0-i\pi,$$
which takes value in $\C/(\cl(\gamma)\Z + 2\pi i\Z)$. We can similarly define the shear between two points in $N^1(\sqrt{\gamma})$ for a cuff $\gamma$.

\subsubsection{Hamster wheels}

The hamster wheels serve as important building blocks to overcome the difficulties that appear in cusped hyperbolic $3$-manifolds. In this paper, we use them only indirectly, by invoking results from \cite{KW21}. We therefore recall only their basic definition.

For a sufficiently large integer $R$, the $R$-perfect hamster wheel $H_R$ is the unique planar hyperbolic surface with geodesic boundary satisfying the following:
\begin{enumerate}
    \item it has $R+2$ geodesic boundary components, each of length $2R$;
    \item it admits an isometric $\Z/R\Z$ action, which fixes two of the boundaries setwise and cyclically permutes the other $R$ boundaries.
\end{enumerate}
Informally, a good hamster wheel can be viewed as a map from $H_R$ to $M$ whose lift to the universal covers is a quasi-isometric embedding with small distortion. In particular, such a map is $\pi_1$-injective. We refer to \cite{KW21} for the precise definition and further details.

We also use, indirectly, the notion of \emph{well-matchedness} for good pants and hamster wheels; see \cite[Section 2.10]{KW21}.

\section{The rose surface}\label{rose-section}

In this section, we introduce the new building block, the rose, to construct surface subgroups. Then we present how to build a \emph{good rose} in cusped hyperbolic manifolds, and state an estimate for the geometry of good roses.

\subsection{The perfect rose}

\begin{figure}[htbp]  
\begin{center}
\begin{overpic}[width=0.3\linewidth]{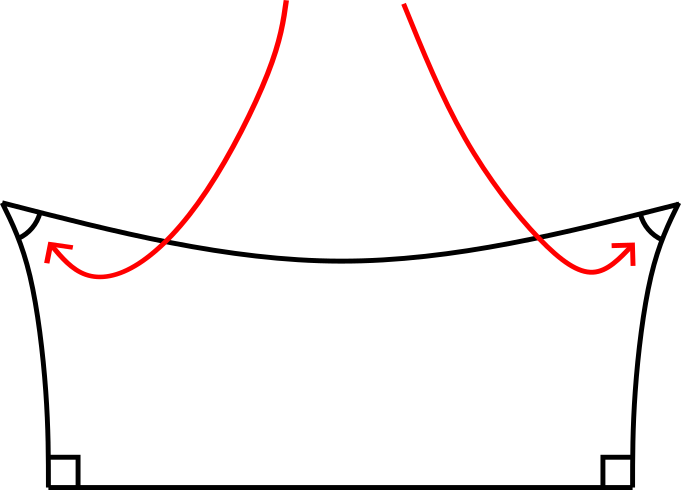}
  \put(37,39){\text{summit}}
  \put(42,3){\text{base}}
  \put(26,77){\textcolor{red}{summit angles}}
  \put(-8,18){\text{leg}}
  \put(97,18){\text{leg}}
\end{overpic}
\caption{A hyperbolic Saccheri quadrilateral.}  
\label{fig:Saccheri}
\end{center}
\end{figure}

Fix a positive integer $g$, which will be determined later. For sufficiently large $L$, let $S_L$ be the hyperbolic Saccheri quadrilateral with summit length $L$ and two summit angles $\pi/4g$; see Figure~\ref{fig:Saccheri}. Then the base and the legs of $S_L$ have length
\begin{equation}
    \widetilde\lambda_L = 2 \arcosh{\left(\sin{\frac{\pi}{4g}}\cosh\frac{L}{2}\right)}, \quad r_L=\artanh\left(\cos\frac{\pi}{4g}\coth\frac{L}{2}\right), 
\end{equation}
respectively.

\begin{figure}[htbp]
\begin{minipage}[t]{0.45\textwidth}
    \centering
    \begin{overpic}[width=0.6\linewidth]{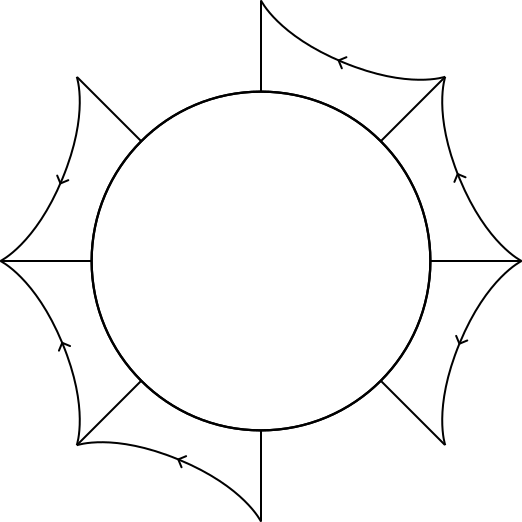}
    \put(91,67){$a_1$}
    \put(65,93){$a_2$}
    \put(2,29){$a_1$}
    \put(29,4){$a_2$}
  
    \put(21,48){$\partial Q_L$}
  
    \put(91,31){$a_{2g}$}
    \put(-2,68){$a_{2g}$}

    \put(63,10){$\cdot$}
    \put(65,10.8){$\cdot$}
    \put(67,11.6){$\cdot$}

    \put(32,84.4){$\cdot$}
    \put(34,85.2){$\cdot$}
    \put(36,86){$\cdot$}
\end{overpic}
    \caption{The perfect rose $Q_L$ as the quotient surface of the annulus.} 
    \label{fig:perfect-rose}
\end{minipage} 
\hspace{0.5cm}
\begin{minipage}[t]{0.45\textwidth}
    \centering
    \begin{overpic}[width=0.6\linewidth]{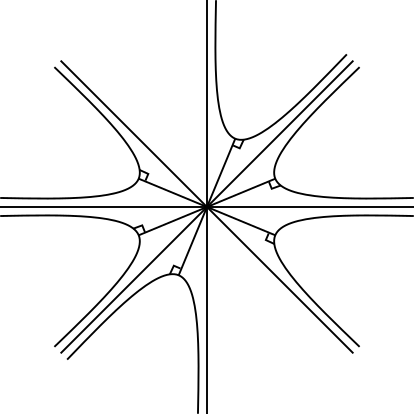}
    \put(81,62){$\widetilde\gamma_1$}
    \put(61,83){$\widetilde\gamma_2$}
    \put(7,36){$\widetilde\gamma_{2g+1}$}
    \put(25,13){$\widetilde\gamma_{2g+2}$}
  
    \put(44.3,56){$\widetilde\star$}
  
    \put(81,34){$\widetilde\gamma_{4g}$}
    \put(11,62){$\widetilde\gamma_{2g}$}

    \put(63,10){$\cdot$}
    \put(65,10.8){$\cdot$}
    \put(67,11.6){$\cdot$}

    \put(32,84.4){$\cdot$}
    \put(34,85.2){$\cdot$}
    \put(36,86){$\cdot$}
\end{overpic}
    \caption{A local picture in the universal cover of $Q_L$ near a lift of $\star$.}
    \label{fig:fundamental-domain}
\end{minipage} 
\end{figure}

By gluing $4g$ copies of $S_L$ cyclically along their legs, we obtain a hyperbolic annulus with
\begin{enumerate}
    \item one geodesic boundary component of length $4g\widetilde\lambda_L$;
    \item one piecewise geodesic boundary component consisting of $4g$ geodesic segments, each of length $L$;
    \item an interior angle of $\pi/2g$ at each vertex of the piecewise geodesic boundary.
\end{enumerate}
We then identify the $4g$ geodesic segments according to the counterclockwise gluing word
\begin{equation}\label{eq:gluing word}
    a_1 a_2\cdots a_{2g} a_1^{-1} a_2^{-1} \cdots a_{2g}^{-1}.
\end{equation}
This identification sends all vertices of the piecewise geodesic boundary to a single point, denoted by $\star$.
We call the resulting hyperbolic surface the $L$-perfect \emph{rose} (of genus $g$), denoted by $Q_L$. We omit the genus from the notation when it is clear from context, as is usually the case. See Figure~\ref{fig:perfect-rose} for the $L$-perfect rose $Q_L$ and the hyperbolic annulus. Since the angle at each vertex is $\pi/2g$, the total angle at $\star$ is $2\pi$. Thus $Q_L$, which is homeomorphic to $S_{g,1}$, is a smooth hyperbolic surface with a geodesic boundary component. We still denote by $a_1,a_2,\dots, a_{2g}$ their images in $Q_L$, and they freely generate $\pi_1(Q_L,\star)$. Thus, the boundary-oriented geodesic $\partial Q_L$ represents the conjugacy class of $$(a_1 a_2\cdots a_{2g} a_1^{-1} a_2^{-1} \cdots a_{2g}^{-1})^{-1}=a_{2g}\cdots a_1 a_{2g}^{-1}\cdots a_1^{-1}.$$ Moreover, we have
\begin{equation}\label{lambdaL}
    \ell(\partial Q_L) = 4g\widetilde\lambda_L = 8g \arcosh{\left(\sin{\frac{\pi}{4g}}\cosh\frac{L}{2}\right)} := \lambda_L.
\end{equation}
The closed geodesic $\partial Q_L$ is also called the \emph{cuff} of $Q_L$.

According to the gluing word \eqref{eq:gluing word}, it is not hard to see that each $a_i$ is a smooth geodesic passing through $\star$. Let $\ww_i$ be the positive tangent vector of $a_i$ at $\star$; then for $1\leq i \leq 2g-1$, we have
\begin{equation}
    \angle (\ww_i,\ww_{i+1})
    =\angle (\ww_{2g},-\ww_{1})
    =\pi-\frac{\pi}{2g}.
\end{equation}
Therefore, in the tangent plane $T_\star Q_L$, the vectors occur counterclockwise in the following cyclic order:
\begin{equation}
    \ww_1, -\ww_2, \ww_3, \dots, \ww_{2g-1}, -\ww_{2g}, -\ww_1, \ww_2,-\ww_3, \dots, -\ww_{2g-1}, \ww_{2g}.
\end{equation}
Among the vectors above, the angle between any two consecutive vectors is $\pi/2g$.

Since $Q_L$ is obtained by gluing $4g$ copies of the Saccheri quadrilateral $S_L$ along legs, the legs then become geodesic segments from $\star$ to $\partial Q_L$ that are orthogonal to $\partial Q_L$. We call these segments the \emph{legs} of $Q_L$, denoted by $\zeta_j$, $1\leq j\leq 4g$. The foot of $\zeta_i$ on $\partial Q_L$ is called the \emph{leg foot} on $\partial Q_L$.
We arbitrarily distinguish a leg foot of $Q_L$ and call it the \emph{formal foot} of $Q_L$.

Now we consider the orthogeodesics from $\partial Q_L$ to itself in $Q_L$. We lift $Q_L$ to its universal cover, denoted by $\widetilde{Q}$, in $\H^2$, and let $\widetilde{\star}$ be a fixed lift of $\star$. Then there are $4g$ lifts $\widetilde\gamma_j$, $1\leq j\leq 4g$, of $\partial Q_L$ whose distances to $\widetilde{\star}$ are equal to the leg length $r_L$ of $S_L$. 
Let $\mathcal B(\widetilde\star)$
be the union of these $4g$ lifts of $\partial Q_L$ joined to $\widetilde\star$ by lifts of the legs of $Q_L$.
See Figure~\ref{fig:fundamental-domain} for the picture described above. Note that any two distinct $\widetilde\gamma_j$ have an orthogeodesic in $\widetilde{Q}$, whose length is
\begin{equation}
    d_k = 2\arcosh\left( \sin\frac{k\pi}{4g} \cosh r_L \right),
\end{equation}
for some $1\leq k\leq 2g$, and this orthogeodesic is a lift of an orthogeodesic from $\partial Q_L$ to itself in $Q_L$. 
Besides, when $2\leq k\leq 2g$, $d_k$ is uniformly bounded above and bounded away from $0$. We now fix a constant $D_0=D_0(g)$ throughout the paper such that
\begin{equation}\label{eq:D_0}
    \frac{1}{D_0} < d_k < D_0,
\end{equation}
for any $2\leq k\leq 2g$ and $L$ sufficiently large.

\begin{figure}[htbp]  
\begin{center}
\begin{overpic}[width=0.3\linewidth]{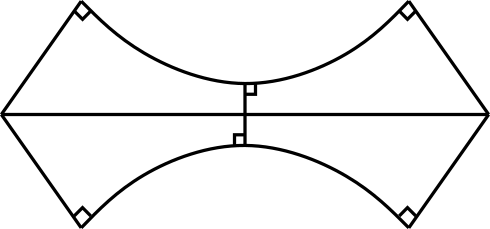}
  \put(46.8,10){$\eta_i$}
  \put(79,28){$a_i$}
\end{overpic}
\caption{A short orthogeodesic $\eta_i$ of $Q_L$ across $a_i$.}  
\label{fig:short-orthogeodesic}
\end{center}
\end{figure}

For $k=1$, there are $2g$ orthogeodesics of length $d_1$, and we call them the \emph{short orthogeodesics} of $Q_L$. They can also be obtained by gluing the two Saccheri quadrilaterals with summits labeled as $a_i$ and then taking the orthogeodesic between the two bases in the resulting hexagon. See Figure~\ref{fig:short-orthogeodesic}.
We denote by $\eta_i$ the short orthogeodesic crossing $a_i$. Then we have 
\begin{equation}
    \ell (\eta_i) = d_1 = 2\arcosh\left( \sin\frac{\pi}{4g} \cosh r_L \right) = 2\artanh\left(\frac{\cot(\pi/4g)}{\sinh(L/2)}\right),
\end{equation}
and therefore $\ell(\eta_i) \asymp_g e^{-L/2}$.

By the symmetry of the gluing word \eqref{eq:gluing word}, we know that the feet of each $\eta_i$ on $\partial Q_L$ are based at a pair of antipodal points. The feet of the $2g$ short orthogeodesics are called the \emph{short feet} of $Q_L$ on $\partial Q_L$. There are $4g$ short feet in total on $\partial Q_L$.

We denote by $\Foot_L\subset N^1(\partial Q_L)$ the collection of short and leg feet of $Q_L$. We remark that the short feet and leg feet occur alternately along $\partial Q_L$. Moreover, the basepoint of any short foot on $\partial Q_L$ is the midpoint of the basepoints of its two adjacent leg feet, and conversely, the basepoint of any leg foot on $\partial Q_L$ is the midpoint of the basepoints of its two adjacent short feet.

\subsection{The good rose}

Given a hyperbolic $3$-manifold $M$, a good rose in $M$ is essentially the image of a perfect rose under a nearly isometric immersion. We now present the more specific definition.

Let $f\from Q_L\to M$ be a continuous map. We regard $f(\star)$ as the basepoint of $M$, and we still denote it by $\star$. The map $f$ (or its image) is said to be an \emph{$(L,\delta)$-good rose}, if there exists an oriented tangent plane $P\subset T_\star M$ at $\star$ and  unit tangent vectors $\u_j\in P$ for $1\leq j\leq 4g$ satisfying that
\begin{enumerate}
    \item $\u_1,\u_2,\dots,\u_{4g}$ occur in this cyclic order counterclockwise, and
    \item for $1\leq j\leq 4g$, $\angle (\u_j,\u_{j+1} ) = \pi/2g$, with the convention $\u_{4g+1}=\u_1$ (thus $\u_i=-\u_{i+2g}$ holds),
\end{enumerate}
such that for $1\leq j\leq 2g$, each $f_{*}(a_j)$ is represented by the carrier of a $\partial$-framed segment ${\frs}_j$ in ${M}$ with the following properties:
\begin{enumerate}
    \item the initial and terminal points of ${\frs}_j$ are both $\star$,

    \item the vectors $\vec{t}_{ini}({\frs}_j), \vec{t}_{ter}({\frs}_j), \vec{n}_{ini}({\frs}_j), \vec{n}_{ter}({\frs}_j)$ are $\de$-close to $(-1)^{j-1}\u_{j},(-1)^{j-1}\u_{j},\n,\n$, respectively,
    
    \item $|\cl({\frs}_j)-L|<\delta$.
\end{enumerate}
Here $\n$ is the positively oriented unit normal vector to $P$.

Recall that $\pi_1(Q_L,\star)=\langle a_1,a_2,\dots, a_{2g}\rangle$ is a free group of rank $2g$. We prove the following lemma.

\begin{lemma}\label{lem:immersed rose}
    For any $\pi/8g>\delta>0$, there exists $L_0>0$ such that for any $L>L_0$, the following holds. Suppose that $f\from Q_L\to M$ is an $(L,\delta)$-good rose, then $f$ is $\pi_1$-injective, and every nonidentity element of the image of $f_*$ is loxodromic.
\end{lemma}

\begin{proof}
    Suppose that $\frs_i$'s are the $\partial$-framed segments associated to the good rose $f$.
    It is sufficient to show that for any nontrivial cyclically reduced word $w\in \pi_1(Q_L)$, its image $f_*(w)$ is loxodromic.
    
    Let $\theta = \pi-\frac{\pi}{2g}+2\de$ and
    \begin{equation}
        w = a_{i_1}^{\de_1} a_{i_2}^{\de_2} \cdots a_{i_t}^{\de_t}
    \end{equation}
    be nontrivial and cyclically reduced, where $t\in\Z_+$ and $\de_k\in\{\pm1\}$.
    Then for any positive integer $n$, $w^n$ is nontrivial and cyclically reduced, and we have
    \begin{equation}
        (f_*(w))^n = (f_*(a_{i_1})^{\de_1} f_*(a_{i_2})^{\de_2} \cdots f_*(a_{i_t})^{\de_t})^n.
    \end{equation}
    Let $\mathfrak{b}_{k}$ represent $(f_*(a_{i_k}))^{\de_k}$, which is either $\frs_{i_k}$ or $\overline{\frs_{i_k}}$.
    By our construction of good roses and the fact that $w$ is cyclically reduced, we know that the bending angle between $\mathfrak{b}_{k}$ and $\mathfrak{b}_{k+1}$ is at most $\theta$, for all $1\leq k\leq t$, with the convention $\mathfrak{b}_{t+1}=\mathfrak{b}_1$. Thus $$\mathfrak{b}_{1},\mathfrak{b}_{2},\dots,\mathfrak{b}_{t},\mathfrak{b}_{1},\mathfrak{b}_{2},\dots,\mathfrak{b}_{t},\dots,\mathfrak{b}_{{1}},\mathfrak{b}_{2},\dots,\mathfrak{b}_{t}$$ form an $(L',\theta)$-tame $2\de$-consecutive chain of oriented $\partial$-framed segments, where $L'=L/2-\delta$.

    Given a lift $\star_1$ of $\star$ in $\H^3$, the universal cover of $M$,
    we lift the above chain
    inductively to $\H^3$ to obtain an $(L',\theta)$-tame $2\de$-consecutive chain of oriented $\partial$-framed segments in $\H^3$, 
    and let $\star_1$, $\star_2$ be the initial and terminal points of the lifted concatenation in $\H^3$.
    Then by Lemma~\ref{lem:ineff} with $\ep=\de$, when $L$ is sufficiently large, we have
    \begin{equation}
        d_\H(\star_1,\star_2) 
        > nt(L-\delta) - (nt-1) I(\theta) - (nt-1)\de > nt (L-2\de-I(\theta)),
    \end{equation}
    for any positive integer $n$.
    Therefore $f_*(w)$ is a loxodromic element.
\end{proof}

Lemma~\ref{lem:immersed rose} implies that when $L$ is sufficiently large, any $(L,\delta)$-good rose is a $\pi_1$-injective immersed subsurface in $M$. Now we introduce the geometry of good roses. Suppose $Q$ is an $(L,\delta)$-good rose given by $f\from Q_L\to M$. We then homotope $f$ so that $f(\partial Q_L)$ is a closed geodesic, denoted by $\partial Q$, in $M$, which is called the \emph{cuff} of $Q$. We estimate the complex length of $\partial Q$.

\begin{lemma}\label{lem:rose-boundary}
    There exists a constant $C$ depending only on $g$, such that for any $\pi/8g>\delta>0$, there exists $L_0>0$ such that for any $L>L_0$, we have 
    \begin{equation}
        |\cl(\partial Q)-\lambda_L|<C\delta.
    \end{equation}
\end{lemma}

\begin{proof}
    We still denote by $\frs_i$ the $\partial$-framed segments associated to the good rose $f$ and let $\theta = \pi-\frac{\pi}{2g}+2\de$. 
    Suppose $c = a_{2g}\cdots a_1 a_{2g}^{-1}\cdots a_1^{-1} \in \pi_1(Q_L)$.
    Then the closed geodesic $\partial Q$ is represented by the reduced cyclic concatenation
    \begin{equation}
        [\frs_{2g}\cdots \frs_1 \overline{\frs}_{2g} \cdots \overline{\frs}_{1}].
    \end{equation}
    By the same argument in the proof of Lemma~\ref{lem:immersed rose}, we know that $\frs_{2g},\cdots, \frs_1, \overline{\frs}_{2g}, \cdots, \overline{\frs}_{1}$ form an $(L',\theta)$-tame $2\de$-consecutive cycle of oriented $\partial$-framed segments, where $L'=L/2-\delta$. 
    
    Denote by $\theta_k$, $1\leq k\leq 4g$, the bending angles of this cycle. Then by construction, we know that each $\theta_k$ is $2\de$-close to $\pi-\frac{\pi}{2g}$. By equation~\eqref{eq:angle-ineff} and $0<\delta<\pi/8g$, we have 
    \begin{equation}\label{eq:angle-diff}
        |I\left(\theta_k\right)-I(\pi-\frac{\pi}{2g})|<C_g\delta,
    \end{equation}
    for some constant $C_g$ depending only on $g$.
    Since $|\cl(\frs_i)-L|<\delta$, we have
    \begin{equation}\label{eq:length-phase}
        |l(\frs_i)-L|, |\phi(\frs_i)| < \delta.
    \end{equation}
    Thus by equations~\eqref{eq:angle-diff}, \eqref{eq:length-phase}, and Lemma~\ref{lem:ineff} with $\ep=\de$, we have
    \begin{equation}\label{eq:boundary-length}
        |l([\frs_{2g}\cdots \frs_1 \overline{\frs}_{2g} \cdots \overline{\frs}_{1}]) - 4gL + 4g\cdot I(\pi-\frac{\pi}{2g})| <  4g\de + 4g C_g\de + 4g\de = (8g+4gC_g)\de
    \end{equation}
    and
    \begin{equation}\label{eq:bounadry-phase}
        |\phi([\frs_{2g}\cdots \frs_1 \overline{\frs}_{2g} \cdots \overline{\frs}_{1}])| < 4g\delta + 8g\de + 4g\delta = 16g\de.
    \end{equation}
    Note that
    \begin{equation}\label{eq:hyperb-trig}
        \widetilde{\lambda}_L = 2 \arcosh{\left(\sin{\frac{\pi}{4g}}\cosh\frac{L}{2}\right)} = L - 2 \ln \left(\sec\left(\frac{\pi}{2}-\frac{\pi}{4g}\right)\right) + O(e^{-L}).
    \end{equation}
    Hence by equations~\eqref{eq:boundary-length}, \eqref{eq:bounadry-phase}, and \eqref{eq:hyperb-trig}, when $L$ is sufficiently large, we have
    \begin{equation}
        |\cl(\partial Q)-\lambda_L| = |\cl(\partial Q)-4g\widetilde{\lambda}_L| < (24g + 4gC_g + 1)\delta. \qedhere
    \end{equation}
\end{proof}

Now, we define the feet of $Q$. For the short orthogeodesic $\eta_i$ of $Q_L$, we homotope $f$ so that each $f(\eta_i)$ is an orthogeodesic $\beta_i$ from $\partial Q$ to itself. We call the $\beta_i$ the \emph{short orthogeodesics} of $Q$.
The two feet of $\beta_i$ are called the \emph{short feet} of $Q$ on $\partial Q$.
We can furthermore homotope $f$ so that each $f(\zeta_j)$ is a geodesic segment $\xi_j$ from the basepoint $\star$ of $Q$ to $\partial Q$ that is orthogonal to $\partial Q$. The $\xi_j$ are called the \emph{legs} of $Q$, and the foot of each $\xi_j$ on $\partial Q$ is called a \emph{leg foot} on $\partial Q$. 
Then there is a unique leg foot of $Q$ that corresponds to the formal foot of $Q_L$, and we call this leg foot the \emph{formal foot} of $Q$.
For a good rose $Q$, we denote by $\Foot(Q)\subset N^1(\partial Q)$ the collection of short and leg feet of $Q$.

To understand the geometry of $Q$, we define the complex length of its short orthogeodesic. Let $\beta$ be a short orthogeodesic of $Q$, which is the image of an orthogeodesic $\eta$ of $Q_L$ under $f$. Let $w_0,w_1$ be two unit tangent vectors of $\partial Q_L$ based at the two endpoints of $\eta$, respectively, so that the lifts of $w_0$ and $w_1$ point to the same side of the complete geodesic containing the lift of $\eta$ after lifting to the universal cover of $Q_L$. Let $w'_0,w'_1$ be the corresponding unit tangent vectors of $\partial Q$ based at the two endpoints of $\beta$, respectively; then $w'_0,w'_1$ are normal to $\beta$. We define $\cl(\beta)$ to be the complex distance between $w'_0$ and $w'_1$; here $\beta$ is oriented so that $\cl(\beta)$ has positive real part. Note that $\cl(\beta)$ does not depend on the choices of $w_0$ and $w_1$.

We now state the following important theorem concerning the geometry of good roses. Since its proof is rather technical, we defer it to Appendix A.

\begin{theorem}\label{thm:good-rose-geometry}
    There exist constants $C,\de_0>0$ depending only on $g$ such that for any $0<\de<\de_0$, there exists $L_0=L_0(g,\de)>0$ such that the following holds for all $L>L_0$. Suppose $Q$ is an $(L,\de)$-good rose. Then we have:
    \begin{enumerate}[label=(\arabic*)]
        \item For each short orthogeodesic $\eta$ of $Q_L$ and its corresponding orthogeodesic $\beta$ in $Q$,
        \begin{equation}
            \left| \frac{\cl(\beta)}{\ell(\eta)} -1 \right| < C\de.
        \end{equation}

        \item For any $\widehat  u_1,\widehat  u_2\in\Foot_L$ and their corresponding feet $u_1,u_2\in\Foot(Q)$, we have
        \begin{equation}
            \left|
            d_{\partial Q}(u_1,u_2)
            -
            d_{\partial Q_L}(\widehat  u_1,\widehat  u_2)
        \right|
        <C\delta.
        \end{equation}
        Here $d_{\partial Q_L}(\widehat  u_1,\widehat  u_2)$ denotes the signed distance from $\widehat  u_1$ to $\widehat  u_2$ along $\partial Q_L$, while $d_{\partial Q}(u_1,u_2)$ denotes the corresponding complex distance along $\partial Q$, with compatible lifts chosen according to the marked intervals along $\partial Q_L$.

        \item There is a unique slow and constant turning normal field $v$ on $\partial Q$ passing through the formal foot of $Q$, such that 
        for any $u\in\Foot(Q)$, we have
        \begin{equation}
            \left|u-v(x_{u})\right|<C\de,
        \end{equation}
        where $x_{u}\in\partial Q$ is the basepoint of
        $u$.
    \end{enumerate}
\end{theorem}

\section{Construction of the desired closed surface}\label{sec:construction}

In this section, we construct the desired closed subsurface in $M$ and prove Theorem~\ref{main2}, assuming the injectivity result established in Appendix~\ref{incompressible-section}. 
We first find a good rose in $M$, and then close it up using good pants and good hamster wheels, via nearly geodesic matching. 

\subsection{Immersed rose in \texorpdfstring{$M$}{M}}\label{immersedrose}

In this subsection, we prove the existence of good roses in cusped hyperbolic $3$-manifolds.
We fix the cusped hyperbolic $3$-manifold $M=\H^3/G$, the finite group $F$, and the epimorphism $\rho\from G\to F$. Recall that there is a universal constant $\ep_0$, independent of $G$, such that in the universal cover $\H^3$, the horoballs of distinct cusps consisting of points with injectivity radius at most $\ep_0$ are pairwise disjoint; see, for instance, \cite[Chapter D]{BP92}. 
For a point in $M$, its \emph{height} is defined to be the signed distance to the boundary of these horoballs so that the height is positive if the point lies in the union of these horoballs. 
We define the \emph{height} of a closed geodesic in $M$ as the maximal height among points on the geodesic, and the \emph{height} of a pants or a hamster wheel in $M$ as the maximal height of its cuffs. 

\begin{theorem}\label{thm:immersed rose}
    For any $\delta>0$, if $L>0$ is sufficiently large, there exists a $\pi_1$-injective $(L,\delta)$-good rose $Q$ given by $f\from Q_L\to M$ such that $\rho\circ f_*$ is surjective.
\end{theorem}

\begin{proof}
Denote by $\star$ the basepoint of $M$, which is assumed to have negative height, and we identify $G$ with $\pi_1(M,\star)$. Let $G_1=\ker(\rho)$ and $|F|=n+1$ for some $n\geq0$. Then there exists a finite covering space $\pi\from \widetilde{M}\to M$ such that $\pi_1(\widetilde{M},\star_0)=G_1$, where $\star_0$ is a lift of $\star$ and serves as the basepoint of $\widetilde{M}$. 
Let $F=\{e, x_1,x_2,\dots,x_{n}\}$. Then for each $1\leq i\leq n$, $\star$ has a lift $\star_i\neq \star_0$ such that for any arc $\alpha$ with initial point $\star_0$ and terminal point $\star_i$, we have $\rho([\pi(\alpha)])=x_i$. Moreover, we know that $\star_i$, $0\leq i\leq n$ are all the lifts of $\star$ in $\widetilde{M}$.

Now fix a positive integer $g$ such that $2g> n$, and we construct a good rose of genus $g$ in $M$.
We arbitrarily choose an oriented tangent plane $P\subset T_\star M$ and then unit tangent vectors $\u_j\in P$ for $1\leq j\leq 4g$ such that
\begin{enumerate}
    \item $\u_1,\u_2,\dots,\u_{4g}$ occur in this cyclic order counterclockwise; and
    \item $\angle (\u_j,\u_{j+1}) = \pi/2g$, with the convention $\u_{4g+1}=\u_1$, for any $j$.
\end{enumerate}
Note that for any $1\leq j\leq 2g$, $\u_j=-\u_{j+2g}$ holds.
We denote by $\n\in T_\star M$ the positively oriented unit normal vector to $P$. 
Note that the covering map $\pi$ induces a covering map  $d\pi\from T\widetilde{M} \to TM$.
Then for each $0\leq i\leq n$ and $1\leq j\leq 4g$, denote by $\u_{i,j}$ and $\n_i$ the lift of $\u_{j}$ and $\n$, respectively, based at $\star_i\in\widetilde{M}$.

By \cite[Theorem 4.11]{KW19} and \cite[Theorem 3.2]{KW21}, for any given $\delta>0$, there exists $L_1>L_0>0$ (here, $L_0$ is from Lemma~\ref{lem:immersed rose}), such that for any $L>L_1$, there exist $\partial$-framed segments $\widetilde{\frs}_j$, $1\leq j\leq 2g$, in the covering manifold $\widetilde{M}$, such that:
\begin{enumerate}
    \item The initial point of $\widetilde{\frs}_j$ is $\star_0$, for $1\leq j\leq 2g$.
    
    \item The terminal point of $\widetilde{\frs}_j$ is $\star_j$ if $1\leq j\leq n$, and is $\star_0$ if $n< j\leq 2g$.

    \item The vectors $\vec{t}_{ini}(\widetilde{\frs}_j), \vec{t}_{ter}(\widetilde{\frs}_j), \vec{n}_{ini}(\widetilde{\frs}_j), \vec{n}_{ter}(\widetilde{\frs}_j)$ are $\de$-close to $\mathbf{a},\mathbf{b},\mathbf{c},\mathbf{d}$, respectively, where
    \begin{equation}
        (\mathbf{a},\mathbf{b},\mathbf{c},\mathbf{d}) = 
        \begin{cases}
            (\u_{0,j}, \u_{j,j}, \n_0,\n_j), & 1\leq j\leq n \text{ and } j \text{ is odd}\\
            (\u_{0,j+2g},\u_{j,j+2g},\n_0,\n_j), & 1\leq j\leq n \text{ and } j \text{ is even},\\
            (\u_{0,j},\u_{0,j},\n_0,\n_0), & n< j\leq 2g \text{ and } j \text{ is odd},\\
            (\u_{0,j+2g}, \u_{0,j+2g}, \n_0, \n_0), & n< j\leq 2g \text{ and } j \text{ is even}.
        \end{cases}
    \end{equation}
    
    \item For each $1\leq j\leq 2g$, $|\cl(\widetilde{\frs}_j)-L|<\delta$.

    \item For each $1\leq j\leq 2g$, the height of $\widetilde\frs_j$ is at most $\log L$.
\end{enumerate}
For $1\leq j\leq 2g$, let $$\frs_j = (\pi (s(\widetilde{\frs}_j)), d\pi(\vec{n}_{ini}(\widetilde{\frs}_j)), d\pi(\vec{n}_{ter}(\widetilde{\frs}_j)))$$
be the $\partial$-framed segment in $M$ that is the image of $\widetilde{\frs_j}$ under $\pi$ (and $d\pi$). Here, $s(\widetilde{\frs}_j)$ is the carrier segment of $\widetilde{\frs}_j$. 

Since $\pi(\star_i)=\star$ for any $0\leq i\leq n$, we know that $\star$ is the initial and terminal point of all $\frs_j$. Therefore, we can regard the carrier segment of each $\frs_j$ as an element $g_j\in\pi_1(M,\star)$.
By our construction, we know that for $1\leq i\leq n$,
\begin{equation}\label{eq:rho-image}
    \rho(g_i)=x_i.
\end{equation}

Recall that $\pi_1(Q_L,\star)=\langle a_1,a_2,\dots, a_{2g}\rangle$ is a free group of rank $2g$. By the argument in Lemma~\ref{lem:immersed rose}, we know that there exists a $\pi_1$-injective $(L,\delta)$-good rose $f\from Q_L\to M$ such that $f_*(a_i)=g_i$ for $1\leq i\leq 2g$. By equation~\eqref{eq:rho-image} and the fact that $2g>n$, we have that $\rho\circ f_*$ is surjective.
\end{proof}

\begin{remark}
    With the chosen marking, the oriented boundary represents the conjugacy class of
\begin{equation}
    f_*((a_1 a_2\cdots a_{2g} a_1^{-1} a_2^{-1} \cdots a_{2g}^{-1})^{-1})=g_{2g}\cdots g_1 g_{2g}^{-1}\cdots g_1^{-1}.
\end{equation}
\end{remark}

\subsection{Constructing the closed immersed surface}\label{subsec:construction}

In this subsection, we construct the desired closed immersed subsurface $\Sigma$ in $M$. 

We first recall the following simplified version of \cite[Theorem 3.11]{Sun23}.

\begin{theorem}\label{boundingsurface}
If $\epsilon>0$ is small enough and $R>0$ is a large enough integer, then for any null-homologous oriented $(R,\epsilon)$-multicurve $\Gamma$ in $M$ with height at most $4\log{R}$, there exists a (possibly disconnected) immersed subsurface $f\from S\to M$ bounded by two copies of $\Gamma$ such that the following hold.
\begin{enumerate}
\item $S$ is an oriented nearly geodesic subsurface obtained by gluing $(R,\epsilon)$-good pants and $(R,\epsilon)$-good hamster wheels so that their cuffs are well-matched.
\item If we write $2\Gamma=\gamma_1\cup \cdots \cup \gamma_k$ as a union of its components (thus $k$ is even), then the oriented boundary $\partial S=c_1\cup \cdots \cup c_k$ has $k$ components such that $f$ maps each $c_i$ to $\gamma_i$.
\item Each $c_i$ is contained in a pair of pants $\Pi_i\subset S$ such that $\Pi_i\ne \Pi_j$ for $i\neq j$ and $f|_{\Pi_i}:\Pi_i \to M$ is an $(R,\epsilon)$-good pants in item (1).
\item For any $v_i\in N^1(\sqrt{\gamma_i})$, we can choose $S$ in $M$ such that $|\mathbf{foot}_{\gamma_i}(\Pi_i,c_i)-v_i|<\epsilon$ holds for all $i$.
\end{enumerate}

\end{theorem}

Now for any $\ep>0$, we choose $\delta>0$ so that $C\delta<\ep$, for the constant $C$ in Lemma~\ref{lem:rose-boundary}. We then choose the sufficiently large $L$ in Theorem~\ref{thm:immersed rose} so that $R={\lambda_L}/{2}$ is an integer, where $\lambda_L$ is given in equation~\eqref{lambdaL}.
By Theorem~\ref{thm:immersed rose}, we have a $\pi_1$-injective $(L,\delta)$-good rose $f\from Q_L\to M$, denoted by $Q$ for simplicity, such that $\rho\circ f_*$ is surjective.
Then Lemma~\ref{lem:rose-boundary} implies that $\partial Q$ is an $(R,\epsilon)$-good curve. 

By \cite[Lemma 3.7]{Sun22}, $\partial Q$ lies in the 1-neighborhood of the union $\cup\frs_i$.
In the argument in the proof of Theorem~\ref{thm:immersed rose}, the height of each $\frs_i$ is chosen to be at most $\log{L}$, and hence the height of $\partial Q$ is no greater than $4\log{R}$. Moreover, since $\partial Q$ is the only boundary component of $Q$, $\partial Q$ is null-homologous in $M$.

So $\Gamma=\overline{\partial Q}$, the orientation reversal of $\partial Q$, satisfies the hypotheses of Theorem~\ref{boundingsurface}. By Theorem~\ref{thm:good-rose-geometry} (3), there is a slow and constant turning normal field $v$ on $\partial Q$ that approximately gives the inward tangent direction of $Q$. 
Let $u_0\in N^1(\partial Q)$ be the formal foot of $Q$ and take $v_0\in N^1(\sqrt{\partial Q})$ to be the quotient of $-u_0$.

Then we apply Theorem~\ref{boundingsurface} to construct a nearly geodesic immersed subsurface $S'$ in $M$ whose boundary is $2\partial Q$ and 
the feet on both boundary components are close to $v_0$. Let $\Sigma'$ be obtained from one copy of $S'$ and two copies of $Q$, gluing along the two copies of $\partial Q$. 
Since $\Sigma'$ is possibly disconnected, we let $\Sigma$ be one connected component of $\Sigma'$ that contains at least one copy of $Q$, and let $S=\Sigma\cap S'$ be the part assembled by pants and hamster wheels.
The subsurface $\Sigma$ is nearly geodesic near the copies of $\partial Q$, since $v_0$ points outward with respect to $Q$ and points inward with respect to $S$.
Moreover, by the reasoning as in \cite[Remark 2.12]{Sun25}, we can apply the process in Step IV of \cite[page 558-559]{Sun15} such that any essential arc in $S$
either contains an essential subarc inside a component of $S$ from one boundary component to itself, or intersects at least $e^{R}$ components of $S$.

Note that the surface $\Sigma$ is obtained by gluing $(R,\ep)$-good pants, $(R,\ep)$-good hamster wheels, and one or two $(L,\delta)$-good roses in $M$. Hence there is a hyperbolic surface $\Sigma_0$, obtained by gluing $R$-perfect pants, $R$-perfect hamster wheels, and $L$-perfect roses, such that $\Sigma$ is realized by $f\from \Sigma_0\to M$, and the map $f$ is obtained by gluing the maps associated to the good pants, good hamster wheels, and good roses.
We claim that $f\from\Sigma_0\to M$ is $\pi_1$-injective. More details on the construction of $\Sigma_0$ and the proof of this technical result will be given in Appendix~\ref{incompressible-section}.

\begin{theorem}\label{injective}
The subsurface $f\from\Sigma_0\to M$ is $\pi_1$-injective, and $f_*(\pi_1(\Sigma_0))$ is quasi-Fuchsian.
\end{theorem}

Now we are ready to prove Theorem~\ref{main2}.

\begin{proof}[Proof of Theorem~\ref{main2}]
We choose the basepoint $\star_\Sigma\in \Sigma_0$ to be the basepoint $\star$ of one copy of the perfect rose $Q_L\subset \Sigma_0$, and take the basepoint $\star_M$ to be $f(\star_{\Sigma})$. 
We then let $H=f_*(\pi_1(\Sigma_0))$, and Theorem~\ref{injective} implies that $H$ is a quasi-Fuchsian surface subgroup of $\pi_1(M)$. 

For the good rose $f\from (Q_L,\star)\to (M,\star_M)$, we know that the composition $\pi_1(Q_L,\star)\to \pi_1(M,\star_M)=G\xrightarrow{\rho}F$ is surjective. Since $Q_L$ is a subsurface of $\Sigma_0$, we have a sequence of injections $$\pi_1(Q_L)\to \pi_1(\Sigma_0)\xrightarrow{f_*}f_*(\pi_1(\Sigma_0))=H\to \pi_1(M)=G.$$ Then $\rho(H)=F$ holds since $\rho$ is already an epimorphism on the subgroup $f_*(\pi_1(Q_L))<H$. 
\end{proof}

\section{Proof of Theorem \ref{main}}\label{provemain}

We prove that Theorem \ref{main2} implies Theorem \ref{main} in this section. 
We begin with the following proposition, whose proof follows the argument of \cite[Proposition~7.1]{KRS}, where a related statement is proved under arithmetic hypotheses. The method goes back to \cite{McReynolds}. We include the proof for the reader's convenience.

\begin{proposition}\label{torsionfree}
Let $G$ be a finite-volume hyperbolic $3$-orbifold group, and let $K<G$ be a finitely generated free subgroup. Then $K$ is contained in a torsion-free finite-index subgroup $G_0<G$.
\end{proposition}

\begin{proof}

Since $G$ is a finite-volume hyperbolic $3$-orbifold group, it has only finitely many conjugacy classes of nonidentity torsion elements. Let $g_1,\cdots, g_n$ be representatives of these conjugacy classes.

Since $K$ is finitely generated, any finite-index subgroup of $K$ is still finitely generated. Since $G$ is LERF (\cite{Agol}), all these finite-index subgroups are separable in $G$. So the inclusion $i: K\to G$ induces an injective homomorphism $\widehat{i}: \widehat{K}\to \widehat{G}$ on their profinite completions. We still use $\widehat{K}$ to denote its image in $\widehat{G}$. 
Since $K$ is a finitely generated free group, $\widehat{K}$ is a finitely generated free profinite group, so it is torsion-free. 

For any $i=1,\cdots, n$, let $$C_i=\{\hat{g}\cdot g_i\cdot \hat{g}^{-1}\ |\ \hat{g}\in \widehat{G}\}.$$ Since $\widehat{G}$ is compact, $C_i$ is a compact subset of $\widehat{G}$ consisting of torsion elements. Since $\widehat{K}$ is torsion-free, we have $\widehat{K}\cap C_i=\varnothing$. Since $\widehat{G}$ is compact and Hausdorff, $\widehat{K}$ and $C_i$ have disjoint open neighborhoods $V_i\supset \widehat{K}$ and $W_i\supset C_i$. 

By the definition of the profinite topology on $\widehat{G}$, there exists an open (and finite-index) normal subgroup $U_i\lhd \widehat{G}$ such that $\widehat{K}\cdot U_i\subset V_i$ and $C_i\cdot U_i\subset W_i$. Then we have $(\widehat{K}\cdot U_i)\cap (C_i\cdot U_i)=\varnothing$. 

Let $\widehat{q}_i:\widehat{G}\to \widehat{G}/U_i$ be the natural quotient homomorphism. We take $$G_i=G\cap \widehat{q}_i^{-1}(\widehat{q}_i(\widehat{K})),$$ and clearly $K<G_i$ holds.
Since $\widehat{G}/U_i$ is a finite group, $\widehat{q}_i^{-1}(\widehat{q}_i(\widehat{K}))$ is a finite-index subgroup of $\widehat{G}$, and $G_i$ is a finite-index subgroup of $G$. Since $(\widehat{K}\cdot U_i)\cap (C_i\cdot U_i)=\varnothing$, we have $\widehat{q}_i(\widehat{K})\cap \widehat{q}_i(C_i)=\varnothing$. So $\widehat{q}_i^{-1}(\widehat{q}_i(\widehat{K}))$ contains no element in $C_i$, thus $g_i$ does not conjugate into $G_i$ (within $G$).

Let $G_0$ be the intersection of all $G_1,\cdots, G_n$. Then it is a finite-index subgroup of $G$. Since $K<G_i$ for each $i$, we know that $K<G_0$. For any $i=1,\cdots, n$, since $g_i$ does not conjugate into $G_i$ (within $G$), it does not conjugate into $G_0$ either (within $G$). Since any nonidentity torsion element in $G$ is conjugate to one of $g_1,\cdots,g_n$, $G_0$ is torsion-free. The subgroup $G_0$ satisfies all the desired conditions of this proposition.
\end{proof}

Now we are ready to prove Theorem \ref{main}.
\begin{proof}[Proof of Theorem \ref{main}]
{\bf Step I.} We first prove Theorem \ref{main} in the case of $A=F$. Let $O=\H^3/G$ be the corresponding 3-orbifold.

By Selberg's lemma, $G$ contains a finite-index torsion-free subgroup $G'$, and we can assume that $G'<\text{ker}(\rho)$ by passing to a further finite-index subgroup. Then $M'=\mathbb{H}^3/G'$ is a finite-volume hyperbolic $3$-manifold. 
Choose a regular basepoint of $O$ outside the cusp horoballs, if any, and
choose lifts of this point to $M'$ representing all elements of $F$. Repeating the connection construction in the proof of Theorem \ref{thm:immersed rose}, we obtain a free subgroup $K<G$ such that $\rho(K)=F$.

By Proposition~\ref{torsionfree}, $K$ lies in a torsion-free finite-index subgroup $G_0<G$. Since $\rho(K)=F$, we have $\rho(G_0)=F$.
Now we apply Theorem \ref{main2} to the epimorphism $\rho|_{G_0}:G_0\to F$. There exists a surface subgroup $H<G_0$ such that $\rho|_{G_0}(H)=F$. Then $H$ is also a surface subgroup of $G$ and $\rho(H)=F$ holds.

{\bf Step II.} Now we prove Theorem \ref{main} for the general case.

For the epimorphism $\rho:G\to F$ and subgroup $A<F$, we take $G'=\rho^{-1}(A)$. Then we have an epimorphism $\rho|_{G'}:G'\to A$. Since $F$ is a finite group, $G'$ is a finite-index subgroup of $G$ and is also a finite covolume Kleinian group. By Step I, there exists a surface subgroup $H<G'$ such that $\rho|_{G'}(H)=A$. Then $H$ is also a surface subgroup of $G$, and $\rho(H)=A$ holds.
\end{proof}

We make the following remark on the abundance of surface subgroups in Theorem \ref{main}.
\begin{remark}\label{abundance}
In the setting of Theorem \ref{main}, we can actually find infinitely many commensurability classes of surface subgroups $H<G$, such that $\rho(H)=A$.

At first, by the above proof of Theorem \ref{main}, we can reduce to the case that $A=F$ and $G=\pi_1(M)$ is a $3$-manifold group. The surface subgroup $H<G$ that surjects onto $F$ is realized by a nearly geodesic subsurface $\Sigma_0\to M$. Let $c_0$ be a decomposition curve of $\Sigma_0$, such that any arc that connects $c_0$ and any rose component of $\Sigma_0$ must intersect at least ${e^R}/{2}$ components of $\Sigma_0$. The existence of $c_0$ is guaranteed by the construction before Theorem \ref{injective}. 

Let $\gamma$ be the image of $c_0$, then it is a good curve. We take another nearly geodesic subsurface obtained by the Kahn--Wright construction $\Sigma_1\to M$ that contains a decomposition curve $c_1$, such that $c_1$ is mapped to $\gamma$, and $\Sigma_1$ intersects with $\Sigma_0$ along $\gamma$ at an angle of almost ${\pi}/{2}$. Let $R\subset \Sigma_0$ be one rose component. We construct a finite cover $\Sigma_0'\to \Sigma_0$ by taking a finite quotient of $\pi_1(\Sigma_0)\to \pi_1(\Sigma_0/R)\to Q$, in which $R$ lifts. Let $c_0'$ be one elevation of $c_0$ in $\Sigma_0'$, and let $d$ be the covering degree of $c_0'\to c_0$. Then we take a finite cover $\Sigma_1'\to \Sigma_1$ such that one elevation $c_1'$ of $c_1$ also covers $c_1$ by degree $d$.

Then we apply the cut-and-paste argument in \cite[Section 4.1]{KM12a} to $c_0'\subset \Sigma_0'$ and $c_1'\subset \Sigma_1'$, to build a new $\pi_1$-injective immersed subsurface $\Sigma'\to M$. Let $H'<G$ be the surface subgroup corresponding to $\Sigma'$, then it also surjects onto $F$ since it contains the rose $R$. Moreover, the argument in \cite[Section 4]{KM12a} implies that, when we vary the choice of finite covers $\Sigma_0'\to \Sigma_0$ and $\Sigma_1'\to \Sigma_1$, the surface subgroup $H'$ can realize infinitely many commensurability classes of surface subgroups of $G$ that surject onto $F$.

\end{remark}

\section{Application to Kent--Leininger 4-manifolds}\label{sec:KL}

A challenging problem in geometric topology is to understand how homology changes when passing to finite covers. If 
$p:\bar M \rightarrow M$ is a finite cover, one always has a transfer map $\tau : H_k(M; \mathbb Z) \rightarrow H_k(\bar M; \mathbb Z)$,
which satisfies $p_*\circ \tau = \deg(p)\cdot id$. This immediately shows that the transfer map is rationally injective, and hence that the ordinary
Betti numbers can only increase when passing to finite covers. If we have an oriented closed immersed submanifold $N \rightarrow M$ representing some 
homology class $[N]\in H_k(M; \mathbb Z)$, the image of the transfer map $\tau[N]$ naturally breaks up into a sum $\sum [N_i]$, where each $N_i$ 
is an immersed submanifold consisting of a connected lift of $N$. Geometrically, it is interesting to understand
the individual homology classes $[N_i]$.

Let us now specialize to the case where $M$ is a closed hyperbolic $3$-manifold. From the work of Agol--Wise \cite{Agol}, we know that any such $M$ virtually
fibers over $S^1$, and hence has a finite cover $p:\bar M \rightarrow M$ with $\beta_1(\bar M)=\beta_2(\bar M) \neq 0$. Note that the fiber surface 
$i: N \hookrightarrow \bar M$ represents an element of infinite order in $H_2(\bar M)$. If one takes the original $M$ to be a closed hyperbolic 
$3$-manifold which is a rational homology sphere (see for instance the discussion in \cite{FKR23}), then we have an immersed closed subsurface $j:N \hookrightarrow \bar M \to M $ such that $j_*[N]=0 \in H_2(M)$, but where the homology class $i_*[N]$ of a connected lift has infinite order. This example shows that connected lifts/elevations of a surface
can exhibit drastically different homological behavior. As a consequence of our main theorem, we obtain the following:

\begin{elevation}
Let $M$ be a finite volume hyperbolic $3$-manifold, and $p:\bar M\rightarrow M$ a finite cover. Then there exists a $\pi_1$-injective immersed connected 
oriented closed surface $j: N \rightarrow M$, such that there is a unique elevation $i:\bar N \rightarrow \bar M$. As a result, for this surface $N$, if 
$j_*[N]=0\in H_2(M)$ then $i_*[\bar N]=0\in H_2(\bar M)$.
\end{elevation}

\begin{proof}
Let $G=\pi_1(M)$ be the associated finite covolume Kleinian group, and let $\bar G \leq G$ be the finite index subgroup associated to the cover $\bar M$. Then there
exists a further finite index normal subgroup $K=\bigcap g\bar G g^{-1} \triangleleft G$ obtained by intersecting all conjugates of $\bar G$. This gives rise to a 
homomorphism $\rho: G \rightarrow F := G/K$ to a finite group, and a subgroup $A:= \rho (\bar G)$ with the property that $\bar G = \rho^{-1}(A)$. Theorem~\ref{main} then
provides us with a $\pi_1$-injective immersed surface $N \rightarrow M$ whose fundamental group $H:= \pi_1(N)$ has the property that $\rho(H)=F$. 
This implies equality of the indices of subgroups $[H: \bar G \cap H]=[ G: \bar G]$, which establishes the first claim.

For the second claim, we note that in degree two, homology coincides with bordism, so $H_2(M)=\Omega _2^{SO}(M)$. The condition $j_*[N]=0\in H_2(M)$ then implies 
there exists a compact oriented $3$-manifold with boundary $(W, \partial W)$, and a map $\phi: W\rightarrow M$ such that $\phi$ restricts to $\partial W=N$ as $j$. We then consider an elevation $\bar \phi: (\bar W, \partial \bar W) \rightarrow \bar M$ of $\phi$, which will map a finite cover of $(W, \partial W)$ to $\bar M$.  Since $\rho(H)=F$, 
we also have $\rho \big( \phi_\sharp(\pi_1(W))\big)=F$. It follows that the induced finite cover $\bar W$ has the same degree $[G:\bar G]$. Since the restriction of the covering
map to $\partial W$ also has that same degree, we see that $\bar W$ has a single boundary component, which is a finite cover $\bar N$ of N. We conclude that
$i_*[\bar N]=0\in H_2(\bar M)$, as desired.
\end{proof}

Let us now discuss an application of these results to the Kent--Leininger $4$-manifolds, whose construction (see \cite{KL}) we now briefly review. Let $M_8$ denote the figure $8$ 
knot complement, a non-compact finite volume hyperbolic $3$-manifold, and set $\Gamma=\pi_1(M_8)$. In their paper, Kent--Leininger constructed
a {\it type-preserving} embedding $\Delta: \Gamma \rightarrow \Mod(S_{1,3})$ to the pure mapping class group of the thrice-punctured torus. 
We use the embedding $\Delta$ to identify $\Gamma$ with its image subgroup in $\Mod(S_{1,3})$.
The type-preserving condition means that every non-peripheral element of $\Gamma$ maps to a pseudo-Anosov class. Then associated to any quasi-Fuchsian closed subsurface $\Sigma$ in $M_8$, one can form a surface bundle $S_{1,3} \rightarrow E \rightarrow \Sigma$, with holonomy given by $\Delta|_{\pi_1(\Sigma)}$. This is 
a non-compact surface bundle, as the fibers are $S_{1,3}$. The type-preserving condition ensures that the only $\mathbb Z^2$-subgroups of $\pi_1(E)$ are peripheral, i.e. 
are supported in the ends of the $4$-manifold $E$. 

In order to obtain {\it compact} atoroidal surface bundles over surfaces, Kent--Leininger viewed $S_{1,3}$ as a torus with three ordered marked points, and considered
finite branched covers $p:\hat S \rightarrow S_{1,3}$, branched over the marked points. An {\it acceptable} branched cover is one where the local degree at every 
pre-image of a marked point is at least two. When $\hat S$ is an acceptable branched cover,
there is a finite index subgroup $\Mod_p(S_{1,3}) \leq \Mod(S_{1,3})$ which has a type-preserving embedding into $\Mod(\hat S)$, i.e. pseudo-Anosov elements are mapped to 
pseudo-Anosov elements. This finite index subgroup gives rise to a finite index subgroup $\Gamma_p :=\Gamma \cap \Mod_p(S_{1,3})$, which inherits a type-preserving
embedding into $\Mod(\hat S)$. Then for each quasi-Fuchsian closed surface subgroup $\pi_1(\Sigma_g) \leq \Gamma_p$, one obtains a closed $4$-manifold that 
fibers as an $\hat S$-bundle over $\Sigma_g$, with holonomy given by the embedding into $\Mod(\hat S)$.

\begin{corollary}   
For any acceptable branched cover $p: \hat S \rightarrow S_{1,3}$, there are infinitely many surface subgroups $\pi_1(\Sigma) \leq \Gamma_p$ for which the associated
surface bundle $\hat S \rightarrow E\rightarrow \Sigma$ has signature zero.
\end{corollary}

\begin{proof}
The finite index subgroup $\Gamma_p \leq \Gamma$ gives rise to a finite cover $\pi\from M_p \rightarrow M_8$ of the figure $8$ knot complement. Then Corollary \ref{elevation} and Remark \ref{abundance} provide
us with infinitely many $\pi_1$-injective closed, connected quasi-Fuchsian subsurfaces $N\rightarrow M_8$ whose elevation $\Sigma \to M_p$ is also connected. 
Since $M_8$ is a knot complement, $H_2(M_8)=0$, so our corollary also tells us $[\Sigma]=0\in H_2(M_p)=\Omega_2^{SO}(M_p)$. Taking $W$ to be 
a compact oriented 3-manifold $W$, with $\partial W=\Sigma$, equipped with a map $W\to M_p$ extending $\Sigma\to M_p$,
the associated $\hat S$-bundle over $W$ is a $5$-manifold with boundary $E$. This forces $E$ to have signature zero.
\end{proof}

This corollary should be contrasted with the results in Lafont--Miller--Ruffoni \cite{LMR}. In that paper, infinitely many acceptable branched covers $p_i:\hat S_i \rightarrow S_{1,3}$ 
were constructed, with the additional property that the entire $\Mod(S_{1,3})$ had a type-preserving embedding into $\Mod(\hat S_i)$, i.e. without having to pass to a finite
index subgroup. As a result $\Gamma_{p_i} = \Gamma$, and all surface subgroups in $\Gamma$ gave rise to $\hat S_i$-bundles with zero signature. So that result provided
infinitely many (but not all!) fiber surfaces, and arbitrary base surfaces, with signature zero. The result in our corollary instead 
applies to every acceptable branched cover and provides infinitely many
(but not all!) base surfaces for which the associated bundle has signature zero. As such our corollary can be viewed as a complementary result to Lafont--Miller--Ruffoni.

\appendix

\makeatletter
\renewcommand*{\theHsection}{appendix.\Alph{section}}
\renewcommand*{\theHequation}{\theHsection.\arabic{equation}}
\makeatother

\section{The geometry of good roses}\label{appen:good-rose-geometry}

In this section, we prove Theorem~\ref{thm:good-rose-geometry}. We will use the hyperbolic cosine and sine laws for right-angled hexagons with complex edge lengths. Quadrilaterals with two or three right angles may be viewed as degenerate right-angled hexagons.
See, for instance, \cite[Section VI]{Fen89} for further details.

\begin{goodrosegeometry}
    There exist constants $C,\de_0>0$ depending only on $g$ such that for any $0<\de<\de_0$, there exists $L_0=L_0(g,\de)>0$ such that the following holds for all $L>L_0$. Suppose $Q$ is an $(L,\de)$-good rose. Then we have:
    \begin{enumerate}[label=(\arabic*)]
        \item For each short orthogeodesic $\eta$ of $Q_L$ and its corresponding orthogeodesic $\beta$ in $Q$, we have
        \begin{equation}
            \left| \frac{\cl(\beta)}{\ell(\eta)} -1 \right| < C\de.
        \end{equation}

        \item For any $\widehat  u_1,\widehat  u_2\in\Foot_L$ and their corresponding feet $u_1,u_2\in\Foot(Q)$, we have
        \begin{equation}
            \left|
            d_{\partial Q}(u_1,u_2)
            -
            d_{\partial Q_L}(\widehat  u_1,\widehat  u_2)
        \right|
        <C\delta.
        \end{equation}
        Here $d_{\partial Q_L}(\widehat  u_1,\widehat  u_2)$ denotes the signed distance from $\widehat  u_1$ to $\widehat  u_2$ along $\partial Q_L$, while $d_{\partial Q}(u_1,u_2)$ denotes the corresponding complex distance along $\partial Q$, with compatible lifts chosen according to the marked intervals along $\partial Q_L$.

        \item There is a unique slow and constant turning normal field $v$ on $\partial Q$ passing through the formal foot of $Q$, such that 
        for any $u\in\Foot(Q)$, we have
        \begin{equation}
            \left|u-v(x_{u})\right|<C\de,
        \end{equation}
        where $x_{u}\in\partial Q$ is the basepoint of
        $u$.
    \end{enumerate}
\end{goodrosegeometry}

\begin{proof}[Proof of Theorem~\ref{thm:good-rose-geometry}]
Throughout the proof, $C>0$ denotes a constant depending only on $g$,
whose value may change from line to line. 
For simplicity, we let $\alpha=\frac{\pi}{4g}$, $N=4g$, and
\begin{equation}
    b_L=\frac{\widetilde{\lambda}_L}{2}= \arcosh\left(\sin\alpha\cosh\frac L2\right).
\end{equation}
Thus $\lambda_L = \ell(\partial Q_L)=2Nb_L.$
Let $h_L$ be the distance from the summit to the base in the
Saccheri quadrilateral $S_L$; then $h_L=d_1/2$ is half of the length of short orthogeodesics of $Q_L$. Hence, $h_L\asymp_g e^{-L/2}$.

We first derive the local geometry of $Q$ from the data of the defining $\partial$-framed segments. For the perfect rose $Q_L$, we cut $Q_L$ along the $2g$ inner geodesics that represent $a_i\in \pi_1(Q_L,\star)$ for $1\leq i\leq 2g$ and obtain the annular fundamental region, whose piecewise geodesic boundary component consists of $N=4g$ oriented geodesic segments of length $L$. We then let 
\begin{equation}
    \sigma_1,\sigma_2,\dots,\sigma_N
\end{equation}
be the corresponding oriented occurrences of the carriers of the $\partial$-framed segments in the definition of $Q$, following the cyclic order. We lift the piecewise geodesic to $\H^3$, and obtain a periodic chain of infinitely many $\partial$-framed segments. There is an initial frame and a terminal frame at each vertex, and we only keep the initial frame.

Let
\begin{equation}
    X(z)=\begin{pmatrix}
        e^{z/2}&0\\
        0&e^{-z/2}
    \end{pmatrix} \quad\text{and}\quad R=
    \begin{pmatrix}
        \sin\alpha&\cos\alpha\\
        -\cos\alpha&\sin\alpha
    \end{pmatrix}.
\end{equation}
With respect to the above chosen frames, by the construction of $(L,\de)$-good roses, the right action of the transition across the $\sigma_i$ followed by crossing the following vertex has the form $B_i=X(z_i)S_i$,
where
\begin{equation}\label{eq:z-S-control}
    |z_i-L|<C\de \quad\text{and}\quad \|S_i-R\|<C\de.
\end{equation}

Let 
\begin{equation}
    S_i=\begin{pmatrix}
        a_i&b_i\\
        c_i&d_i
    \end{pmatrix}.
\end{equation}
By equation~\eqref{eq:z-S-control}, we have
\begin{equation}\label{eq:rose-S-entries}
    a_i=\sin\alpha+O_g(\de),\quad
    b_i=\cos\alpha+O_g(\de),\quad
    c_i=-\cos\alpha+O_g(\de),\quad
    d_i=\sin\alpha+O_g(\de).
\end{equation}
In particular, after decreasing $\de_0$ if necessary, all $a_i$ and $c_i$ are uniformly bounded away from zero.

We analyze the geometry near $\sigma_1$, and the same argument applies cyclically to other $\sigma_i$.  Let
\begin{equation}
    G=B_1B_2\cdots B_N
    =
    \begin{pmatrix}
        A&B\\
        C&D
    \end{pmatrix}.
\end{equation}
The isometry represented by $G$ is the deck transformation
corresponding to $\partial Q$, and hence its axis is the lift of $\partial Q$ accompanied by the lift of the piecewise geodesic. Let $E_i=e^{z_i/2}$. Since $N=4g$ is fixed, a direct inductive computation gives that
\begin{align}
    A
    &=
    \left(\prod_{i=1}^N E_i a_i\right)
    \bigl(1+O_g(e^{-L})\bigr),
    \label{eq:rose-matrix-A}\\
    \frac BA
    &=
    \frac{b_N}{a_N}
    (1+O_g(e^{-L})),
    \label{eq:rose-matrix-B}\\
    \frac CA
    &=
    e^{-z_1}\frac{c_1}{a_1}
    \bigl(1+O_g(e^{-L})\bigr),
    \label{eq:rose-matrix-C}\\
    \frac DA
    &=O_g(e^{-L}).
    \label{eq:rose-matrix-D}
\end{align}

In particular, for $L$ sufficiently large, equations~\eqref{eq:rose-matrix-A}--\eqref{eq:rose-matrix-D} imply that $G$ is loxodromic, as indicated by Lemma \ref{lem:immersed rose}.
We normalize $\widehat{\sigma}_1$, the lift of $\sigma_1$, so that it lies on the geodesic with endpoints $0,\infty$ in the upper half-space model, with initial point at height $1$, its initial frame given by the normal vector $(1,0)\in\C\times\R$, and its terminal point has height greater than $1$. 
Let $\zeta_-$ and $\zeta_+$ be the endpoints of the axis of $G$. Then $\zeta_-$ and $\zeta_+$ are the roots of
\begin{equation}
    C\zeta^2+(D-A)\zeta-B=0,
\end{equation}
and hence
\begin{equation}\label{eq:rose-fixed-product}
    \zeta_-\zeta_+=-\frac BC,
    \quad
    \zeta_-+\zeta_+=\frac{A-D}{C}.
\end{equation}

Let $m_1$ be the complex distance along $\widehat{\sigma}_1$ from its initial vertex to the foot of its common orthogonal with the axis of $G$. By hyperbolic geometry of the two geodesics $(0,\infty)$ and $(\zeta_-,\zeta_+)$, we have $e^{2m_1}=\zeta_-\zeta_+$, see Figure~\ref{fig:m_1}.
Together with equations~\eqref{eq:rose-fixed-product},
\eqref{eq:rose-matrix-B}, and \eqref{eq:rose-matrix-C}, we obtain
\begin{align*}
    e^{2m_1}
    =e^{z_1}
    \left(-\frac{b_Na_1}{a_Nc_1}\right)\bigl(1+O_g(e^{-L})\bigr)
    =e^{z_1}\bigl(1+O_g(\de)\bigr),
\end{align*}
when $L$ is sufficiently large.
Here the last equality follows from
equation~\eqref{eq:rose-S-entries}.

\begin{figure}[htbp]  
\begin{center}
\begin{overpic}[width=0.4\linewidth]{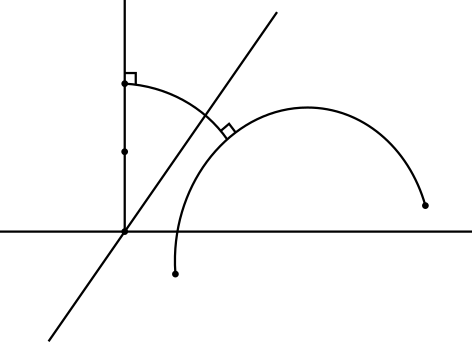}
  \put(38,7){$\zeta_-$}
  \put(93,26){$\zeta_+$}
  \put(17,17){$0$}
  \put(19.5,35){$1$}
  \put(17,47){$m_1$}
  \put(17,70){$\infty$}
  \put(35,55){$h_1$}
\end{overpic}
\caption{The common orthogeodesic between $\widehat{\sigma}_1$ and the axis of $G$.}  
\label{fig:m_1}
\end{center}
\end{figure}

We now claim that 
\begin{equation}\label{eq:rose-midpoint-control}
    \left|m_1-\frac{z_1}{2}\right|
    <C\de,
\end{equation}
and hence, by equation~\eqref{eq:z-S-control}, we have
\begin{equation}\label{eq:rose-midpoint-L}
    \left|m_1-\frac L2\right|
    <C\de.
\end{equation}
Indeed, set $w=e^{2m_1-z_1}$. Then for sufficiently small $\de_0$, we may assume that $|w-1|<1/4$. 
Since $m_1$ and $z_1$ vary continuously with respect to the inputs $z_i$ and $S_i$, 
we have $2m_1-z_1=\Log(w)$, where $\Log(1)=0$.
Hence, equation~\eqref{eq:rose-midpoint-control} holds when $\de_0$ is sufficiently small.

Let $h_1$ denote the complex length of the common orthogeodesic between $\widehat{\sigma}_1$ and the axis of $G$. By the cross-ratio formula, we have
$$\cosh^2 h_1 = \frac{(\zeta_-+\zeta_+)^2}{(\zeta_+-\zeta_-)^2}.$$
Therefore,
$$\tanh^2 h_1=\frac{4\zeta_-\zeta_+}{(\zeta_-+\zeta_+)^2}.$$
Then by equations~\eqref{eq:rose-matrix-C} and \eqref{eq:rose-matrix-D}, we have
\begin{equation}
\begin{aligned}
    \tanh^2 h_1
    &=4\frac{\frac{-B}{A}\cdot\frac{C}{A}}{(1-\frac{D}{A})^2}
    =4(-\cot\alpha+O_g(\de))(-e^{-z_1}(\cot\al+O_g(\de)))(1+O_g(e^{-L}))^2\\
    & = 4e^{-z_1}\cot^2\al(1+O_g(\de))= 4e^{-L}\cot^2\al(1+O_g(\de)).
\end{aligned}
\end{equation}
On the other hand, we have
\begin{equation}
    \tanh h_L=\frac{\cot\alpha}{\sinh(L/2)}=2e^{-L/2}\cot\alpha\bigl(1+O_g(e^{-L})\bigr).
\end{equation}
It follows that
\begin{equation}
    \frac{\tanh^2 h_1}{\tanh^2 h_L}=1+O_g(\de).
\end{equation}
Since $h_1$ has positive real part and $h_L>0$, we have
\begin{equation}\label{eq:rose-tanh-height-relative}
    \frac{\tanh h_1}{\tanh h_L}=1+O_g(\de).
\end{equation}
Since $h_1,h_L=O_g(e^{-L/2})$ and $h_1$ varies continuously with respect to the $z_i$ and $S_i$, we know that when $\de_0$ is sufficiently small,
\begin{equation}\label{eq:rose-height-relative-new}
    \frac{h_1}{h_L}=1+O_g(\de).
\end{equation}

Applying the same argument to each $\sigma_i$, we see that complex parameters $m_i,h_i$ satisfy
\begin{equation}\label{eq:rose-all-local}
    \left|m_i-\frac L2\right|<C\de,
    \quad
    \left|\frac{h_i}{h_L}-1\right|<C\de.
\end{equation}

For $i\in\Z/N\Z$, the geodesic segment $\sigma_i$ has two legs from its endpoints to $\partial Q$. Together with $\sigma_i$ and a piece from $\partial Q$, they form a quadrilateral close to a Saccheri quadrilateral. We call it the $i$-th Saccheri-type piece. We let $x_i,y_i$ be the complex summit angles at the initial and terminal point of $\sigma_i$ of the $i$-th Saccheri-type piece. See Figure~\ref{fig:Saccheri-type}.

\begin{figure}[htbp]  
\begin{center}
\begin{overpic}[width=0.4\linewidth]{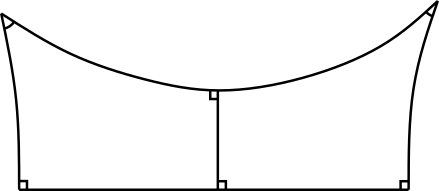}
  \put(5,30){$x_i$}
  \put(90,30){$y_i$}
  \put(24,30){$m_i$}
  \put(42,11){$h_i$}
  \put(60,32){$z_i-m_i$}
  \put(68,3){$b_i^+$}
  \put(25,3){$b_i^-$}
\end{overpic}
\caption{The $i$-th Saccheri-type piece.}  
\label{fig:Saccheri-type}
\end{center}
\end{figure}

Note that the common orthogonal between $\sigma_i$ and $\partial Q$ divides the $i$-th Saccheri-type piece into two complex Lambert quadrilaterals. Let $b_i^-$ and $b_i^+$ be the complex distances along $\partial Q$ from the initial leg foot to the foot of the common orthogonal on $\partial Q$ and from the foot of the common orthogonal to the terminal leg foot.
Then $x_i,y_i$ and $b_i^\pm$ also vary continuously with respect to the $z_i$ and $S_i$.
For $x_i$ and $y_i$, we have
\begin{equation}
    \cot x_i=\tanh h_i\,\sinh m_i,
    \quad
    \cot y_i=\tanh h_i\,\sinh(z_i-m_i).
\end{equation}
By equations~\eqref{eq:z-S-control}, \eqref{eq:rose-all-local}, $\cot\alpha=\tanh h_L\,\sinh(L/2)$, and the same continuity argument, we have 
\begin{equation}\label{eq:rose-angle-control-new}
    |x_i-\alpha|, |y_i-\alpha|<C\de.
\end{equation}
For $b_i^\pm$, we have
\begin{equation}
    \cosh b_i^-=\sin x_i\cosh m_i,
    \quad
    \cosh b_i^+=\sin y_i\cosh(z_i-m_i).
\end{equation}
By $\cosh b_L=\sin\alpha\cosh(L/2)$, equations~\eqref{eq:z-S-control},
\eqref{eq:rose-all-local}, and
\eqref{eq:rose-angle-control-new}, the same continuity argument implies
\begin{equation}\label{eq:rose-boundary-half-control}
    |b_i^\pm-b_L|<C\de.
\end{equation}

We are now ready to prove the theorem.
To see Part (1),
let $\eta$ be a short orthogeodesic of $Q_L$ and
let $\beta$ be the corresponding orthogeodesic in $Q$. The two long sides paired to form $\eta$ are identified by one of the generators of the rose. After making this identification in the universal cover, let $h_-$ and $h_+$ be the complex distances from the resulting common long geodesic $\sigma$ to the two relevant lifts of $\partial Q$. By equation~\eqref{eq:rose-height-relative-new},
\begin{equation}\label{eq:rose-hpm}
    h_\pm=h_L\bigl(1+O_g(\de)\bigr).
\end{equation}
Let $s$ be the shear along $\sigma$ between the feet of these two orthogonals as shown in Figure~\ref{fig:two-Saccheri}. The basepoints of both feet lie within $C\de$ of the midpoint of the corresponding long side by
equation~\eqref{eq:rose-midpoint-control}. Besides, at the two endpoints of $\sigma$, the tangent vectors to $\sigma$ and the two legs are nearly coplanar; hence
\begin{equation}\label{eq:rose-s-control-new}
    |s|<C\de.
\end{equation}

\begin{figure}
\begin{center}
\begin{overpic}[width=0.5\linewidth]{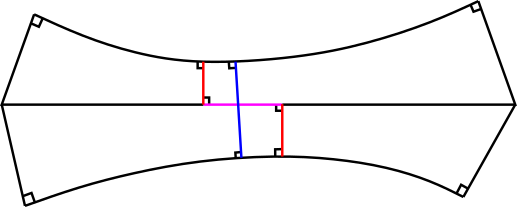}
  \put(32.6,23){\textcolor{red}{$h_+$}}
  \put(56,13){\textcolor{red}{$h_-$}}
  \put(49,21.2){\textcolor[HTML]{FF00FF}{$s$}}
  \put(41,13){\textcolor{blue}{$D$}}
  \put(80,22){$\sigma$}
  \put(48,4.5){$X$}
\end{overpic}
\caption{Two Saccheri-type pieces are glued along the common summit $\sigma$.}  
\label{fig:two-Saccheri}
\end{center}
\end{figure}

Let $D=\cl(\beta)$. By hyperbolic cosine law for right-angled hexagons, we have
\begin{equation}\label{eq:rose-short-hexagon-new}
    \cosh D = \cosh h_-\cosh h_+ +\sinh h_-\sinh h_+ \cosh s.
\end{equation}
Equivalently,
\begin{equation}
    \sinh^2\frac D2
    = \sinh^2\frac{h_-+h_+}{2} + \frac12 \sinh h_-\sinh h_+(\cosh s-1).
\end{equation}
By equations~\eqref{eq:rose-hpm} and \eqref{eq:rose-s-control-new}, the second term is $O_g(h_L^2\de^2)$, while
\begin{equation}
    \left|\sinh\frac{h_-+h_+}{2}\right|\asymp_g h_L.
\end{equation}
Note again that $D$ is continuous with respect to the $z_i$ and $S_i$. Then it follows that
\begin{equation}\label{eq:D}
    D= h_-+h_+ + O_g(h_L\de^2)
    = 2h_L\bigl(1+O_g(\de)\bigr).
\end{equation}
Since $\ell(\eta)=2h_L$, we have
\begin{equation}
    \left|\frac{\cl(\beta)}{\ell(\eta)}-1 \right|
    <C\de,
\end{equation}
which proves Part (1).

To see Part (2), we need the following estimate for the feet. Let $q$ be the foot on one of the boundary lifts of the auxiliary common orthogonal in the proof of Part (1) of length $h_-$, and let $u$ be the actual foot of $\beta$ on the same boundary lift. 
Let $X=d_{\partial Q}(q,u)$, which is continuous in the $z_i$ and $S_i$, see Figure~\ref{fig:two-Saccheri}.
By hyperbolic sine law for right-angled hexagons, we have
\begin{equation}
    \sinh X = \frac{\sinh h_+ \sinh s}{\sinh D}.
\end{equation}
Then by equations~\eqref{eq:rose-hpm}, \eqref{eq:rose-s-control-new}, and \eqref{eq:D}, we have $|\sinh X| = O_g(\de)$.
Therefore, the continuity argument implies that 
\begin{equation}\label{eq:rose-short-foot-control}
    |X|<C\de.
\end{equation}

We recall that the basepoints of the elements of $\Foot_L$ alternate between the short feet and the leg feet and divide $\partial Q_L$ into $2N=8g$ intervals, each of length $b_L$. 
For $Q$, equation~\eqref{eq:rose-boundary-half-control} controls the distance from every leg foot to the corresponding auxiliary midpoint foot, while equation~\eqref{eq:rose-short-foot-control} shows that the latter is within $C\de$ of the actual short foot. Thus every complex distance along $\partial Q$ between two consecutive elements of $\Foot(Q)$ differs from $b_L$ by at most $C\de$.
Since there are only $8g$ such intervals, this proves Part (2).

To see Part (3), we recall that by Lemma \ref{lem:rose-boundary}, $$\left|\cl(\partial Q)-\lambda_L\right|<C\de.$$
Let $\cl(\partial Q)=B+i\theta$. Then $|\theta|<C\de$.

Let $\widehat u_0\in\Foot_L$ be the formal foot of $Q_L$, and let $u_0\in\Foot(Q)$ be the formal foot of $Q$ with basepoint $x_0$. Then there is a unique slow and constant turning vector field $v$ on $\partial Q$ such that $v(x_0)=u_0$.

Suppose $u\in\Foot(Q)$ is an arbitrary foot. Let
$\widehat  u\in\Foot_L$ be the corresponding object of $Q_L$ and let $x$ be the basepoint of $u$. 
Denote by $s_0=d_{\partial Q_L}(\widehat u_0, \widehat u)$, and let $s=d_{\partial Q}(u_0,u)$ be the compatible lift of the complex distance. By Part (2), we have $|s-s_0|<C\de.$
Since $s_0$ is real, $|\im s|<C\de.$

Let $t=\re s$ be the signed distance along $\partial Q$ from $x_0$ to $x$. Along the corresponding interval on $\partial Q$, the constant turning field $v$ rotates by $\frac{\theta}{B}t$ between $x_0$ and $x$. Since $|t|\leq B$ and $|\theta|<C\de$, we have
\begin{equation}
    \left|\frac{\theta}{B}t\right|\leq |\theta|
    <C\de.
\end{equation}
Together with $|\im s|<C\de$, it implies
\begin{equation}
    \left|u-v(x)\right|<C\de.
\end{equation}
Since $u$ was arbitrary, this proves Part (3). When $\de_0$ is sufficiently small, we have $$|\theta|<C\de<\pi,$$ so $v$ is indeed the slow and constant turning vector field.
\end{proof}

\section{Incompressibility of the resulting surface}\label{incompressible-section}

In this section, we apply the machinery developed in \cite[Appendix A]{KW21} to show that the resulting surface assembled from good roses, good pants, and good hamster wheels is $\pi_1$-injective. We first recall several important concepts and results. See \cite[Appendix A.1]{KW21} for the motivation.

\subsection{Assemblies and good assemblies}

We now extend the notion of an assembly to include roses by the following definitions, which generalize the corresponding concepts in \cite{KW21}.

For a compact hyperbolic surface with geodesic boundary $Q_0$ and a map $f\from Q_0\to M$, we call it a \emph{component of type} $Q_0$, if $f$ is $\pi_1$-injective and the image of every nontrivial loop under $f$ is homotopic to a closed geodesic. For any component $f$, we can find a representative $\hat{f}$ freely homotopic to $f$ such that $\hat{f}(\gamma)$ is a geodesic for each component $\gamma$ of $\partial Q_0$. 
In this section, we will always assume components are oriented and satisfy this property. 

Suppose $\mathcal F = \{f_i\from Q_i\to M\}$ is a multiset of oriented components, and there is a fixed-point-free involution $\tau$ on the multiset $\cup_i \partial Q_i$ such that for each $\alpha\in \cup_i \partial Q_i$, $\al$ and $\tau(\al)$ are mapped to the same geodesic in $M$ with opposite orientations.
Then we call $\mathcal{A}=(\mathcal{F},\tau)$ an \emph{assembly}. By identifying $\al$ and $\tau(\al)$ via the identification of their associated geodesics in $M$, we obtain a closed oriented topological surface $S_\mathcal{A}=(\cup Q_i)/\tau$, by joining the $Q_i$. 
We let $f_\mathcal{A}\colon S_\mathcal{A}\to M$ be the joining of $f_i$'s via $\tau$.

We consider only pants, hamster wheels, and roses. For $R,\ep>0$, an $(R,\ep)$-\emph{good component} is a pair of $(R,\ep)$-good pants, an $(R,\ep)$-good hamster wheel, or an $(L,\de)$-good rose, where $L,\de>0$ are determined by $\lambda_L=2R$ and $C\de=\ep$ for the constant $C$ from Lemma~\ref{lem:rose-boundary}. Thus the boundary of an $(L,\de)$-good rose is an $(R,\ep)$-good curve.

Given two good components and a curve $\gamma$ that occurs in both of their boundaries, we say that they are \emph{well-glued} if either of the following holds:
\begin{enumerate}
    \item Neither component is a rose. In this case, they are well-matched, as defined in \cite{KW21}.
    \item One component is a good rose, and the other is a good pants. 
    Let $u_0\in N^1(\sqrt{\gamma})$ be the image of the rose’s formal foot, and let $u_1\in N^1(\sqrt{\gamma})$ be the foot of the pants. Then we require that $|u_0-u_1-i\pi|<\ep$.
\end{enumerate}

An assembly $\mathcal{A}$ is $(R,\ep)$-\emph{good}, if the following hold:
\begin{enumerate}
    \item Every component is $(R,\ep)$-good.

    \item For every component $Q$ and a cuff $\gamma$, $Q$ is matched with a well-glued component along $\gamma$.

    \item Let $\mathcal A_1$ be the subassembly consisting of all good pants and hamster wheels in $\mathcal A$. Then every 
    essential arc in $S_{\mathcal A_1}$ that does not contain an essential subarc inside a component of $\mathcal A_1$ from one boundary component to itself, must intersect at least $e^{R}$ components of $\mathcal A_1$.
\end{enumerate}

\subsection{Frames and distortion}

We then recall some notations for the frames.
If $v\in T^1(\H^2)\subset T^1(\H^3)$, then we can complete $v$ to a 2-frame $w,v$ for $\H^2$, and then complete $w,v$ to a 3-frame $q,w,v$ for $\H^3$. This frame is called the associated 3-frame for $v$.
In a $3$-manifold $M$, suppose that $\gamma$ is an oriented geodesic in $M$. Then any unit normal vector $v$ to $\gamma$ determines a unique 2-frame $w,v$ in $M$, where $w$ is the unit positively-oriented tangent vector to $\gamma$. We then complete it to a 3-frame $q,w,v$ in $M$, which is called the associated 3-frame for $v$ with respect to $\gamma$. 
We denote by $\mathcal{F}(\gamma)$ the set of all such frames. If $\gamma$ is unoriented, then we let $\mathcal{F}(\gamma)=\mathcal{F}(\gamma^+)\cup\mathcal{F}(\gamma^-)$, where $\gamma^+$ and $\gamma^-$ are two oppositely oriented versions of $\gamma$. 
For a perfect component $\widehat Q$, we let $\partial^\mathcal{F} \widehat{Q}$ be the set of associated $3$-frames for inward-pointing unit normal vectors based at points in $\partial \widehat{Q}$. For a \emph{good} component $Q$, we let $\partial^\mathcal{F} Q = \cup_{\gamma\in\partial Q}\mathcal{F}(\gamma)$, where $\gamma\in\partial Q$ is oriented.

We then introduce the distance and distortion in $\mathcal{F}\H^3$.
We fix a left-invariant Riemannian metric $d$ on $\Isom^+(\H^3)$, and for $g\in\Isom^+(\H^3)$, let $d(g)=d(1,g)$. For two frames $u,v\in\F\H^3$, we denote by $u\to v \in\Isom^+(\H^3)$ the unique right action determined by $u\cdot (u\to v) = v$. For a subspace $X\subset \F\H^3$ and a map $\phi\from X\to\ \F\H^3$, the map $\phi$ is said to be $\ep$\emph{-distorted up to distance} $D$, if for any $u,v\in X$ and $d(u,v)<D$, one has
\begin{equation}
    d(u\to v, \phi(u)\to \phi(v))<\ep.
\end{equation}
For two maps $\phi_1,\phi_2\from X\to \F\H^3$, we say that $\phi_1$ and $\phi_2$ are $\ep$\emph{-related}, if for any $x\in X$,
\begin{equation}
    d(\phi_1(x)\to\phi_2(x))<\ep.
\end{equation}
We then have the following facts.

\begin{lemma}[{\cite[Lemma~A.8]{KW21}}]
For all $D$, there exists $D'$ such that if $U\in\SL(2,\C)$ and $\Vert U\Vert<D$, then $d(U)<D'$. This is also true if $\Vert\cdot\Vert$ and $d(\cdot)$ are interchanged.
\end{lemma}

\begin{lemma}[{\cite[Lemma~A.9]{KW21}}]
For all $D,\ep$ there exists $\de$ such that if $U,V\in\SL(2,\C)$, $\Vert U\Vert,\Vert V\Vert<D$ and $\Vert U-V\Vert<\de$, then $d(U,V)<\ep$.
\end{lemma}

\begin{lemma}[{\cite[Lemma~A.10]{KW21}}]\label{lem:error-sum}
For all $D,k$ there exists $C=C(D,k)>0$ such that when $u_0,\cdots,u_k$, $v_0 \cdots ,v_k \in \mathcal{F}\H^3$, and
\begin{equation}
    d(u_i\to u_{i+1},v_i\to v_{i+1})<\de,
\end{equation}
and
\begin{equation}
    d(u_i\to u_{i+1})<D,
\end{equation}
then
\begin{equation}
    d(u_0\to u_k,v_0\to v_k)<C\de.
\end{equation}
\end{lemma}

\begin{remark}
    Lemma~\ref{lem:error-sum} is a quantitative strengthening of \cite[Lemma~A.10]{KW21}. It follows immediately from the uniform Lipschitz continuity of right multiplication on bounded subsets of $\Isom^+(\H^3)$, and we leave the details to the reader.
\end{remark}

Now we briefly introduce some formal manipulations in the proof. Suppose that we have a collection of good components $P_i$ (pants, hamster wheels, and roses) glued together along boundaries, such that each cuff of one $P_i$ is geometrically identified with some cuff of another $P_j$, so that we obtain an assembly $\mathcal{A}$. We then construct a perfect model $\mathcal{\widehat{A}}$ for $\mathcal{A}$ that has a perfect component $\widehat{P}_i$ for each $P_i$ in $\mathcal{A}$ along with a map $h_i\colon\widehat{P}_i\to M$ sending $\partial\widehat{P}_i$ to $\partial P_i$ up to homotopy through such maps. We glue the $\widehat{P}_i$'s together, through the isometric identifications of the cuff of the $\widehat{P}_i$, so we obtain a geometric surface $S_{\widehat{\mathcal{A}}}$ and a homotopy class of maps $h\colon S_{\widehat{\mathcal{A}}}\to M$ which maps each $\widehat{P}_i$ to $P_i$ in $M$. Then we lift this map to $\widetilde{h}\colon\H^2\to\H^3$, where the homotopy class of $h$ determines the relationship of each lift of $\widehat{P}_i$ and the corresponding lift of $P_i$, and the same for boundary geodesics of perfect and good components.
After constructing maps $e\colon\partial^\mathcal{F}\widehat{P}_i\to\partial^\mathcal{F}P_i$ that send frames over each boundary geodesic of $\widehat{P}_i$ to frames over the corresponding boundary geodesic of $P_i$, we can use $\widetilde{h}$ to get a canonical lift $\widetilde{e}$ of $e$ to $\partial^\mathcal{F}(\widetilde{\widehat{\mathcal{A}}})$. We then say $e\colon\partial^\mathcal{F}\widehat{\mathcal{A}}\to\partial^\mathcal{F}\mathcal{A}$ is $\ep$-distorted up to distance $D$ if and only if $\widetilde{e}$ does.

\subsection{Compliant maps}

For a good rose (pants, or hamster wheel) $Q$, we let $\widehat{Q}$ be its perfect model. We say that $\phi\from \partial^\F \widehat{Q} \to \partial^\F Q$ is $D,\ep$\emph{-compliant} if the following properties hold.
\begin{enumerate}[label=(\arabic*)]
    \item The map from $\partial \widehat{Q}$ to $\partial Q$ induced by $\phi$ is the restriction (to $\partial \widehat{Q}$) of an orientation-preserving homeomorphism from $\widehat{Q}$ to $Q$.

    \item The induced map is affine on each component of $\partial\widehat{Q}$, and maps each component of $\partial^\F\widehat{Q}$ to the frames determined by a slow and constant turning normal field on $\partial Q$.

    \item The map $\phi$ sends each formal foot of $\widehat{Q}$ to the corresponding formal foot of $Q$.

    \item The map $\phi$ is $\ep$-distorted up to distance $D$.
\end{enumerate}

Before we state the compliance result for good roses, we first study the local geometry of the universal covers of perfect roses.
For a lift
$\widetilde\star$ of $\star$, recall that
$\mathcal B(\widetilde\star)$
is the union of $4g$ lifts of $\partial Q_L$ in the universal cover $\widetilde{Q}_L$ that are joined to $\widetilde\star$ by lifts of the legs of $Q_L$.

\begin{lemma}\label{lem:local-picture}
For every $D>0$, there exist $T_D,L_D>0$ such that, for $L>L_D$, the following holds. Let $\gamma_0,\gamma_1$ be two boundary components of $\widetilde{Q}_L$, and let $w_i\in\partial^{\F}\gamma_i$ such that $d(w_0,w_1)<D.$
Then either $\gamma_0=\gamma_1$, or there exists a lift $\widetilde\star$ of $\star$ such that $\gamma_0,\gamma_1\in\mathcal B(\widetilde\star)$.

In the latter case, let $k\in\{1,\ldots,2g\}$ and
$$\gamma^{(0)}=\gamma_0,\gamma^{(1)},\ldots,\gamma^{(k)}=\gamma_1$$
be the shorter cyclic sequence from $\gamma_0$ to $\gamma_1$ in $\mathcal B(\widetilde\star)$.
Then the consecutive lifts $\gamma^{(i)},\gamma^{(i+1)}$ are joined by lifts of short orthogeodesics. Let $\widehat q_i$ denote the inward-pointing leg-foot frame on $\gamma^{(i)}$ whose leg ends at $\widetilde\star$. Then for $0\leq i<k$,
\begin{equation}\label{eq:qi}
 d(\widehat q_i,\widehat q_{i+1})<T_D.
\end{equation}
Moreover, if $k\geq2$, then
\begin{equation}\label{eq:endpoint}
 d(w_0,\widehat q_0), d(\widehat q_k,w_1)<T_D.
\end{equation}
\end{lemma}

\begin{proof}
    Let $p\colon \widetilde Q_L\longrightarrow Q_L$ be the universal covering map, and set
    $$\mathcal T_L:=p^{-1}(a_1\cup\cdots\cup a_{2g}).$$
    Thus $\mathcal T_L$ is a geodesic tree in $\H^2$ whose vertices are the lifts of $\star$ and whose edges have length $L$. At every vertex, the angle between two distinct incident rays is at least $\pi/2g$. Each lifted Saccheri quadrilateral has one summit as an edge of $\mathcal T_L$ and one base contained in a lift of $\partial Q_L$.

    Since the leg length $r_L$ is uniformly bounded, fix $r>0$ such that $r_L<r$ for all sufficiently large $L$. 
    Consider the summit and the base of a lifted Saccheri quadrilateral. By basic hyperbolic geometry, we know that 
    the Hausdorff distance between the summit and the base is at most $r$.
    Consequently, every point of every lift of $\partial Q_L$ lies within distance $r$ of some edge of $\mathcal T_L$.

    Let $x_i$ be the basepoint of $w_i$. Then there exists $D_1=D_1(D)$ such that when $d(w_0,w_1)<D$, we have $d_{\H^2}(x_0,x_1)<D_1$. Choose an edge $e_i\subset \T_L$ and $y_i\in e_i$ such that $d_{\H^2}(x_i,y_i)<r$. Thus 
    \begin{equation}\label{eq:y_i-1}
        d_{\H^2}(y_0,y_1)<D_1+2r.
    \end{equation}

    On the other hand, we claim that there exists $C_g>0$ such that for sufficiently large $L$, the distance between two points in any two non-incident edges of $\mathcal T_L$ is at least $L-C_g$. Indeed, suppose $e_0,e_1$ are two non-incident edges and $\gamma$ the shortest edge path joining $e_0$ and $e_1$ in $\T_L$. Let $\gamma_1=\gamma\cup e_0\cup e_1$; then $\gamma_1$ is a consecutive chain of geodesic segments of length $L$ and contains at least $3$ edges. Note that the bending angle at each vertex of $\gamma_1$ is at most $\pi-\pi/(2g)$. By Lemma~\ref{lem:ineff}, there exists $C'_g>0$ such that the distance between the two endpoints of $\gamma_1$ is at least $nL-(n-1)C'_g$, where $n\geq 3$ is the number of edges in $\gamma_1$. Therefore for any $y_i\in e_i$, $i=0,1$, we have 
    \begin{equation}\label{eq:y_i-2}
        d_{\H^2}(y_0,y_1)\geq nL-(n-1)C'_g - 2L = (n-2)L - (n-1) C'_g > L - C_g,
    \end{equation}
    for some constant $C_g>0$.

    By equations~\eqref{eq:y_i-1} and \eqref{eq:y_i-2}, we know $e_0$ and $e_1$ are either equal or incident to a common vertex $\widetilde\star$ of $\mathcal T_L$.
    If $e_0=e_1$ and $\gamma_0\neq\gamma_1$, then $\gamma_0$ and $\gamma_1$ are the two boundary lifts adjacent to the two sides of that edge, so they are both joined by legs to either endpoint of $e_0$.
    If $e_0\neq e_1$, then $e_0$ and $e_1$ have a common endpoint, and the lifted Saccheri quadrilaterals containing $x_0$ and $x_1$ both have a leg ending at the common endpoint. Thus, unless $\gamma_0=\gamma_1$, we have $\gamma_0,\gamma_1\in\mathcal B(\widetilde\star)$ for some $\widetilde{\star}$. This proves the first assertion.

    Now suppose $\gamma_0\neq\gamma_1$ and $\gamma_0,\gamma_1\in \mathcal B(\widetilde\star)$.
    By the cyclic geometry around $\widetilde\star$, consecutive lifts are joined by short orthogeodesics, and the shorter cyclic sequence has length $k\leq2g$.
    For consecutive $\gamma^{(i)},\gamma^{(i+1)}$, the basepoints of the leg-foot frames $\widehat q_i,\widehat q_{i+1}$, denoted by $x_i,x_{i+1}$ respectively, are both at distance $r_L$ from $\widetilde\star$, and the angle between the two legs is $\pi/(2g)$. By the hyperbolic cosine law,
    \begin{equation}
        \cosh d_{\mathbb H^2}(x_i,x_{i+1})=\cosh^2 r_L-\sinh^2 r_L\cos\frac{\pi}{2g}.
    \end{equation}
    Since $r_L$ is uniformly bounded, enlarging $T_D$ if necessary gives \eqref{eq:qi}.

    It remains to prove \eqref{eq:endpoint}. Assume $k\ge2$, and let $p_0\in\gamma_0$ and $p_k\in\gamma_1$ be the basepoints of the feet of their common orthogeodesic $\beta$. The length of $\beta$ is $d_k$, with $\frac1{D_0}<d_k<D_0$, where $D_0$ is the fixed constant given by equation~\eqref{eq:D_0}.
    Let $s_0$ be the signed distance along $\gamma_0$ from $p_0$ to the basepoint $x_0$ of $w_0$. The Lambert quadrilateral determined by $x_0$, $p_0$, and $\beta$ gives
    \begin{equation}
        \sinh d_{\mathbb H^2}(x_0,\gamma_1)=\sinh d_k\cosh |s_0|.
    \end{equation}
    Since $d_{\mathbb H^2}(x_0,\gamma_1)\le d_{\mathbb H^2}(x_0,x_1)<D_1$, we obtain
    \begin{equation}\label{eq:s_0}
        \cosh |s_0| \leq \frac{\sinh D_1}{\sinh d_k} \leq \frac{\sinh D_1}{\sinh(1/D_0)}.
    \end{equation}
    Thus $|s_0|$ is bounded by a constant depending only on $D$ and $g$.

    On the other hand, let $q_0\in\gamma_0$ and $q_k\in\gamma_1$ be the basepoints of the feet $\widehat{q}_0$ and $\widehat{q}_k$. Then $p_0,q_0,\widetilde{\star},q_k$, and $p_k$ form a pentagon with four right angles, with the remaining angle $k\pi/2g$. The two sides incident to $\widetilde{\star}$ have length $r_L$. Let $t_k=d_{\gamma_0}(q_0,p_0)$. By hyperbolic trigonometry, we have
    \begin{equation}
        \tanh |t_k| = \frac{\cot\frac{k\pi}{4g}}{\sinh r_L}.
    \end{equation}
    When $L$ is sufficiently large, we have
    \begin{equation}
        r_L=\artanh \left(\cos\frac{\pi}{4g} \coth\frac{L}{2}\right) > \artanh\left(\cos\frac{\pi}{4g}\right).
    \end{equation}
    Therefore $\sinh r_L > \cot (\pi/(4g))$. By $k\geq 2$, we have
    \begin{equation}
        \tanh |t_k|\leq \frac{\cot (2\pi/4g)}{\cot (\pi/4g)}<1.
    \end{equation}
    Hence, we have a uniform upper bound for $|t_k|$:
    \begin{equation}
        |t_k|\leq \artanh\left(\frac{\cot (\pi/2g)}{\cot (\pi/4g)}\right).
    \end{equation}
    Together with equation~\eqref{eq:s_0}, we obtain a uniform bound, depending only on $D$ and $g$, for the distance between $x_0$ and $q_0$. Increasing $T_D$ then gives $d(w_0,\widehat{q}_0)<T_D$. The same argument also works for $w_1$ and $\widehat{q}_k$, and therefore, we complete the proof. 
\end{proof}

Let $\delta_0$ be the constant from Theorem~\ref{thm:good-rose-geometry}. We then prove the following theorem.

\begin{theorem}\label{thm:rose-compliance}
    For all $D>0$, there exists $K>0$ such that for $\de_0>\de>0$, there exists $L_0>0$ such that for $L>L_0$, the following holds. Let $Q$ be an $(L,\de)$-good rose and $\widehat{Q}$ its perfect model, then there is a $D,K\de$-compliant map $\phi\from\partial^\F\widehat{Q}\to\partial^\F Q$.
\end{theorem}

\begin{proof}
    Throughout the proof, constants denoted by $C$ depend only on $g$, whose value may change from line to line, while constants such as $C_D$ may also depend on $D$. We take $L_0$ sufficiently large whenever necessary.

    We recall that $\lambda_L=\ell(\partial Q_L)$ is the length of the boundary of $Q_L=\widehat{Q}$. Let $\cl(\partial Q) = B+i\theta$ be the complex length of $\partial Q$. Then by Lemma~\ref{lem:rose-boundary}, we have
    \begin{equation}\label{eq:cuff}
        |B-\lambda_L|+|\theta|<C\de.
    \end{equation}
    Let $\kappa=\frac{B+i\theta}{\lambda_L}$. By $\lambda_L=4gL+O_g(1)$, we have
    \begin{equation}\label{eq:kappa}
        |\kappa-1|<\frac{C\de}{L}.
    \end{equation}

    Let $\widehat{u}_0$ and $u_0$ be the formal foot of $\widehat{Q}$ and $Q$, respectively. By Theorem~\ref{thm:good-rose-geometry}(3), there is a unique slow and constant turning normal field $v$ on $\partial Q$ that passes through the formal foot $u_0$. 
    Let $\overline{\phi}\from \partial\widehat{Q}\to\partial Q$ be the unique orientation-preserving affine map that sends the basepoint of $\widehat{u}_0$ to the basepoint of $u_0$.
    Then for an inward-pointing frame $\widehat u(x)\in\partial^{\F}Q_L$ based at $x\in\partial Q_L$, we define
    $\phi(\widehat{u}(x))$ to be the associated $3$-frame determined by $v(\overline{\phi}(x))$. Hence, we have $\phi\from \partial^\F \widehat{Q} \to \partial^\F Q$.
    It is straightforward to see that $\phi$ satisfies Conditions (1)-(3) of $D, K\de$-compliance. 

    Let $\widehat u\in\Foot_L$ and let $u\in\Foot(Q)$ be its corresponding foot of $Q$. We claim that
    \begin{equation}\label{eq:feet-close}
        d\bigl(\phi(\widehat u),u\bigr)<C\delta.
    \end{equation}
    Indeed, let $\widehat s=d_{\partial Q_L}(\widehat u_0,\widehat u)\in [0,\lambda_L)$ and $s=d_{\partial Q}(u_0,u)$ be the corresponding complex distance. By Theorem~\ref{thm:good-rose-geometry}(2), we have $|\widehat s - s|<C\de$.
    On the other hand, the affine map $\phi$ sends the displacement $\widehat{s}$ to $\kappa \widehat{s}$, and hence equation~\eqref{eq:kappa} implies $|\kappa \widehat{s}-\widehat{s}|<C\delta$. Together with Theorem~\ref{thm:good-rose-geometry}(3), the basepoint and normal direction of $\phi(\widehat u)$ differ from those of $u$ by $O_g(\delta)$; equivalently, equation~\eqref{eq:feet-close} holds.

    To prove Condition (4) in the definition of $D, K\de$-compliance, we first prove the following lemma.

    \begin{lemma}\label{lem:adjacent}
        For every $D>0$, there exist $C_D,L_D>0$ such that, for $L>L_D$, the following holds. Let $\gamma_0,\gamma_1$ be two lifts of $\partial\widehat Q$ in the universal cover in $\H^2$. 
        Suppose either that $\gamma_0=\gamma_1$, or that $\gamma_0$ and $\gamma_1$ are joined by a lift of a short orthogeodesic. If $w_i\in\partial^{\F}\gamma_i$ and $d(w_0,w_1)<D$, then
        $$d\bigl(w_0\to w_1, \widetilde\phi(w_0)\to\widetilde\phi(w_1)\bigr)<C_D\delta,$$
        where the $\partial^{\F}\gamma_i$ are the lifts of $\partial^{\F}\widehat{Q}$ on $\gamma_i$ and $\widetilde\phi$ is the canonical lift of $\phi$.
    \end{lemma}

    \begin{proof}[Proof of Lemma~\ref{lem:adjacent}]
        If $\gamma_0=\gamma_1$, the conclusion follows from the affine property of $\phi$ and equation~\eqref{eq:kappa}.

        Otherwise, let $\eta$ be the short orthogeodesic one of whose lifts joins $\gamma_0$ to $\gamma_1$, and $\beta$ the corresponding short orthogeodesic in $Q$. By Theorem~\ref{thm:good-rose-geometry}(1), we have
        \begin{equation}
            \left|\frac{\cl(\beta)}{\ell(\eta)}-1\right|<C\de.
        \end{equation}
        By equation~\eqref{eq:feet-close}, the image of the two feet of $\eta$ under $\tilde{\phi}$ are $C\de$-close to the corresponding feet of $\beta$.
        Then the argument in \cite[Theorem A.16]{KW21} for pants works here, which completes the proof.
    \end{proof}

    We now verify Condition (4). We work in the universal cover and let $w_0,w_1\in\partial^{\F}\widetilde Q_L$ satisfy $d(w_0,w_1)<D$. If they lie on the same boundary lift, or on two lifts joined by a short orthogeodesic, Lemma~\ref{lem:adjacent} gives the desired estimate directly. Otherwise, we apply Lemma~\ref{lem:local-picture}. Then there exists $2g\geq k\geq2$, and we have the chain
    $$w_0,\widehat q_0,\widehat q_1,\ldots,\widehat q_k,w_1.$$ 
    We consider the transitions between any two consecutive frames in this chain. 
    Note that the first and last transitions lie on single boundary lifts and have size bounded by $T_D$, while each intermediate transition joins two boundary lifts connected by a short orthogeodesic and has size bounded by $T_D$. By Lemma~\ref{lem:adjacent}, applied with $D=T_D$, we have an $O_{D,g}(\delta)$ error for every transition. Since there are at most $2g+2$ transitions, by Lemma~\ref{lem:error-sum}, there exists a constant $K=K(D,g)$ such that
    $$d\bigl(w_0\to w_1, \widetilde\phi(w_0)\to\widetilde\phi(w_1)\bigr)<K\delta.$$
    Thus $\phi$ is $K\delta$-distorted up to distance $D$, and Condition (4) holds. Hence $\phi$ is $D,K\delta$-compliant.
\end{proof}

Now suppose that $\mathcal{A}$ is a good assembly and $\widehat{\mathcal{A}}$ is a perfect one. We recall that a map $e\from \partial^\mathcal{F}\widehat{\mathcal{A}}\to \partial^\mathcal{F} \mathcal{A}$ is $M,\ep$-compliant, if it satisfies the following conditions.
\begin{enumerate}
    \item For each component $\widehat{Q}$ of $\widehat{\mathcal{A}}$, there is a corresponding component $Q$ of $\mathcal{A}$ such that $e|_{\partial^\mathcal{F}\widehat{Q}}$ is an $M,\ep$-compliant map from $\partial^\mathcal{F}\widehat{Q}$ to $\partial^\mathcal{F} Q$.
    \item If $\widehat{\gamma}$ is a gluing boundary of $\widehat{\mathcal{A}}$, i.e., $\widehat{\gamma}\in\partial \widehat{Q}_0$ and $\widehat{\gamma}\in\partial \widehat{Q}_1$, and $n_0$ and $n_1$ are unit normal vectors at the same basepoint to $\widehat{\gamma}$ pointing inward to $\widehat{Q}_0$ and $\widehat{Q}_1$ respectively, then
    \begin{equation}
        d(n_0\to n_1,e(n_0)\to e(n_1)) < \ep.
    \end{equation}
\end{enumerate}

By the same argument as in the proof of \cite[Theorem A.17]{KW21}, we have the following theorem.

\begin{theorem}\label{thm:assembly-compliant}
    For all $D$, there exists $K,\ep_0>0$ such that for any $0<\ep<\ep_0$, there exists $R_0$ such that for all $R>R_0$, the following holds: For any $(R,\ep)$-good assembly $\mathcal{A}$, there is a perfect assembly $\widehat{\mathcal{A}}$ and a $D,K\ep$-compliant map $e\from \partial^\mathcal{F}\widehat{\mathcal{A}}\to \partial^\mathcal{F} \mathcal{A}$.
\end{theorem}

\begin{proof}
    The proof mainly follows the argument in the proof of \cite[Theorem A.17]{KW21}, and we only need to indicate how every perfect rose is matched to a perfect pants.

    When a good rose is matched to a good pants in $\mathcal A$, we assemble the corresponding perfect models such that the two formal feet are joined with a vanishing shear. Then Condition~(2) of compliance follows from the $(R,\ep)$-goodness of $\mathcal A$ and Theorem~\ref{thm:rose-compliance}.
\end{proof}

We now prove a local distortion bound for the map constructed above.

\begin{theorem}\label{thm:assembly-distortion}
    For all $D>0$, there exists $C>0$ such that for sufficiently small $\ep>0$, there exists $R_0>0$ such that for all $R>R_0$, the following holds. If $\mathcal{A}$ is an $(R,\ep)$-good assembly, then there is a perfect assembly $\widehat{\mathcal{A}}$ and a map $e\from \partial^\mathcal{F}\widehat{\mathcal{A}}\to \partial^\mathcal{F} \mathcal{A}$ that is $C\ep$-distorted up to distance $D$.
\end{theorem}

\begin{proof}
    We fix $D>0$. Then there exists $D_1=D_1(D)>0$ such that for any $x_0,x_1\in\H^2$ with $\widehat{w}_i$ a $3$-frame based at $x_i$, if $d(\widehat{w}_0,\widehat{w}_1)<D$, then $d_{\H^2}(x_0,x_1)<D_1$. Moreover, there exists $D_2=D_2(D)>0$ such that for any $x_0,x_1\in\H^2$ with $\widehat{w}_i$ a $3$-frame based at $x_i$, if $d_{\H^2}(x_0,x_1)<D_1$, then $d(\widehat{w}_0,\widehat{w}_1)<D_2$.
    Let $D'=\max\{D,D_2\}$, which depends only on $D$ (and $g$).
    
    Let $\mathcal A$ be as given in the statement. By Theorem~\ref{thm:assembly-compliant}, we have a perfect assembly $\widehat{\mathcal{A}}$ and a $D',K\ep$-compliant map $e\from \partial^\mathcal{F}\widehat{\mathcal{A}}\to \partial^\mathcal{F} \mathcal{A}$. Here $K$ depends only on $D$ (and $g$).
    Let $\mathcal A_1\subset \mathcal A$ be the subassembly consisting of all good pants and good hamster wheels in $\mathcal A$, and let $S=S_{\mathcal A_1}$ be the assembled surface for $\mathcal A_1$. 
    Let $\widehat{\mathcal A}_1\subset \widehat{\mathcal A}$ be the subassembly corresponding to $\mathcal A_1$.

    We first show that 
    after passing to the universal cover $\widetilde{\widehat{\mathcal A}}$ of the perfect assembly $\widehat{\mathcal A}$, any geodesic segment of length less than $D_1$ can meet at most one lift of a rose, provided that $R$ is sufficiently large.
    By our construction, the good assembly $\mathcal{A}$ and its corresponding $\widehat{\mathcal{A}}$ have the same combinatorial structure and differ only in the geometry of the corresponding components.
    In particular, we know that Condition~(3) in the definition of a good assembly also holds for $\widehat{\mathcal A}$. That says, every essential arc in the assembled surface $S_{\widehat{\mathcal A}_1}$
    either contains an essential subarc inside a component of $\widehat{\mathcal A}_1$ from one boundary component to itself, or intersects at least $e^{R}$ components of $\widehat{\mathcal A}_1$. By \cite[Remark 2.12]{Sun25}, it follows that every essential arc in $S_{\widehat{\mathcal A}_1}$ has length at least $R$.
    Hence when $R$ is sufficiently large, a geodesic segment of length less than $D_1$ cannot contain an essential arc in $S_{\widehat{\mathcal A}_1}$ as a subsegment. Therefore, such a geodesic segment can meet at most one lift of a rose.

    We now 
    let $\widetilde{e}\from \partial^{\mathcal F}(\widetilde{\widehat A})\to \mathcal F\mathbb H^3$ be the canonical lift of $e$. Let $z_i\in\partial^{\mathcal F}(\widetilde{\widehat{\mathcal A}})$ based at $x_i$ with $d(z_0,z_1)<D$. Then we have $d_{\H^2}(x_0,x_1)<D_1$. Let $\al$ be the geodesic segment connecting $x_0$ and $x_1$. By the preceding argument, when $R$ is sufficiently large, $\al$ intersects at most one lift of a rose.

    We now decompose the transition from $z_0$ to $z_1$ along $\alpha$. Choose frames $$u_0=z_0,u_1,\ldots,u_m=z_1$$ with $m\le 5$, so that every consecutive pair is of one of the following three types:
    \begin{enumerate}
        \item both frames lie on boundary lifts belonging to a portion of $\alpha$ contained entirely in a lift of $\widehat{\mathcal A}_1$;
        \item both frames lie on boundary lifts of the same lifted perfect rose;
        \item the two frames are based at the same point of a gluing curve and point inward to the two components on opposite sides of that curve.
    \end{enumerate}
    The bound $m\le5$ comes from the fact that $\alpha$ meets at most one rose: there are at most two pants/hamster-wheel portions, one rose portion, and two gluing transitions. Since the distance between the basepoints of $u_i$ and $u_{i+1}$ is less than $D_1$, then for every $i$, we have $d(u_i\to u_{i+1})<D_2\leq D'$.

    For every $0\le i<m$, we have
    \begin{equation}
        d\bigl(u_i\to u_{i+1},\widetilde e(u_i)\to\widetilde e(u_{i+1})\bigr) < C_2\epsilon
    \end{equation}
    for some $C_2=C_2(D)$.
    For a pair of type~(1), this follows from the argument in \cite[Theorem A.18]{KW21}. For a pair of type~(2), it follows from rose compliance. For a pair of type~(3), it follows from assembly compliance and Theorem~\ref{thm:assembly-compliant}.

    Since $m\leq 5$, Lemma~\ref{lem:error-sum} implies that there exists a constant $C=C(D)$ such that
    \begin{equation}
        d\bigl(z_0\to z_1,\widetilde e(z_0)\to\widetilde e(z_1)\bigr)<C\epsilon,
    \end{equation}
    which completes the proof.
\end{proof}

\subsection{The proof of \texorpdfstring{$\pi_1$}{π₁}-injectivity}

We are now ready to prove Theorem~\ref{injective}. Recall that $X'\subset X$ is said to be $A$-\emph{dense} in a metric space $X$ if $N_A(X')=X$.

\begin{proof}[Proof of Theorem~\ref{injective}]
    By our construction in Section~\ref{subsec:construction}, we know that the surface $\Sigma$ is the assembled surface for an $(R,\ep)$-good assembly $\mathcal{A}$, and the surface $\Sigma_0$ is the assembled surface for a perfect assembly $\widehat{\mathcal{A}}$, where $\widehat{\mathcal{A}}$ is associated with $\mathcal{A}$ as given in Theorem~\ref{thm:assembly-distortion}.

    Note that for an $L$-perfect rose, its leg length has an upper bound independent of $L$. 
    Let $U=\partial^{\mathcal F}(\widetilde{\widehat{\mathcal{A}}})$ for the universal cover $\widetilde{\widehat{\mathcal{A}}}$ of $\widehat{\mathcal{A}}$; then $U$ is $A$-dense in $\mathcal{F}\H^2$ for some positive constant $A=A(g)$. 
    By Theorem~B.9 and \cite[Theorems~A.19 and~A.23]{KW21}, the induced representation is faithful and quasiconformally conjugate to the Fuchsian representation uniformizing $\Sigma_0$. Thus $f$ is $\pi_1$-injective, and its image group is quasi-Fuchsian.
\end{proof}

\bibliographystyle{alpha}
\bibliography{reference}

@article {FKR23,
    AUTHOR = {Ferrari, L. and Kolpakov, A. and Reid, A. W.},
     TITLE = {Infinitely many arithmetic hyperbolic rational homology
              3-spheres that bound geometrically},
   JOURNAL = {Trans. Amer. Math. Soc.},
  FJOURNAL = {Transactions of the American Mathematical Society},
    VOLUME = {376},
      YEAR = {2023},
    NUMBER = {3},
     PAGES = {1979--1997},
      ISSN = {0002-9947,1088-6850},
   MRCLASS = {57K32 (57R40 57R42)},
  MRNUMBER = {4549697},
MRREVIEWER = {Leone\ Slavich},
       DOI = {10.1090/tran/8816},
       URL = {https://doi.org/10.1090/tran/8816},
}

@article{KL,
  title={Atoroidal surface bundles},
  author={Kent, Autumn E and Leininger, Christopher J},
  journal       = {Ann. of Math. (2)},
  FJOURNAL = {Annals of Mathematics. Second Series},
  note          = {To appear},
  eprint        = {2405.12067},
  archivePrefix = {arXiv},
  primaryClass  = {math.GT},
  url           = {https://arxiv.org/abs/2405.12067}
}

@article{LMR,
  title={On signatures of the atoroidal bundles of {K}ent--{L}eininger},
  author={Lafont, Jean-Fran{\c{c}}ois and Miller, Nicholas and Ruffoni, Lorenzo},
  journal={arXiv preprint arXiv:2410.18029},
  year={2024}
}

@article{KW21,
    AUTHOR = {Kahn, Jeremy and Wright, Alex},
     TITLE = {Nearly {F}uchsian surface subgroups of finite covolume
              {K}leinian groups},
   JOURNAL = {Duke Math. J.},
  FJOURNAL = {Duke Mathematical Journal},
    VOLUME = {170},
      YEAR = {2021},
    NUMBER = {3},
     PAGES = {503--573},
      ISSN = {0012-7094,1547-7398},
   MRCLASS = {57M07 (30C62 30F40)},
  MRNUMBER = {4255043},
MRREVIEWER = {Bruno\ P.\ Zimmermann},
       DOI = {10.1215/00127094-2020-0049},
       URL = {https://doi-org.proxy.libraries.rutgers.edu/10.1215/00127094-2020-0049},
}

@article{KM12b,
    AUTHOR = {Kahn, Jeremy and Markovi\'c, Vladimir},
     TITLE = {Immersing almost geodesic surfaces in a closed hyperbolic
              three manifold},
   JOURNAL = {Ann. of Math. (2)},
  FJOURNAL = {Annals of Mathematics. Second Series},
    VOLUME = {175},
      YEAR = {2012},
    NUMBER = {3},
     PAGES = {1127--1190},
      ISSN = {0003-486X,1939-8980},
   MRCLASS = {57M50 (30F40)},
  MRNUMBER = {2912704},
MRREVIEWER = {James\ W.\ Anderson},
       DOI = {10.4007/annals.2012.175.3.4},
       URL = {https://doi-org.proxy.libraries.rutgers.edu/10.4007/annals.2012.175.3.4},
}

@article{KM12a,
    AUTHOR = {Kahn, Jeremy and Markovi\'c, Vladimir},
     TITLE = {Counting essential surfaces in a closed hyperbolic
              three-manifold},
   JOURNAL = {Geom. Topol.},
  FJOURNAL = {Geometry \& Topology},
    VOLUME = {16},
      YEAR = {2012},
    NUMBER = {1},
     PAGES = {601--624},
      ISSN = {1465-3060,1364-0380},
   MRCLASS = {57M50},
  MRNUMBER = {2916295},
MRREVIEWER = {Colin\ C.\ Adams},
       DOI = {10.2140/gt.2012.16.601},
       URL = {https://doi-org.proxy.libraries.rutgers.edu/10.2140/gt.2012.16.601},
}

@article{KM15,
    AUTHOR = {Kahn, Jeremy and Markovi\'c, Vladimir},
     TITLE = {The good pants homology and the {E}hrenpreis conjecture},
   JOURNAL = {Ann. of Math. (2)},
  FJOURNAL = {Annals of Mathematics. Second Series},
    VOLUME = {182},
      YEAR = {2015},
    NUMBER = {1},
     PAGES = {1--72},
      ISSN = {0003-486X,1939-8980},
   MRCLASS = {57N05 (30C62 30F20 53D25)},
  MRNUMBER = {3374956},
MRREVIEWER = {Bruno\ Martelli},
       DOI = {10.4007/annals.2015.182.1.1},
       URL = {https://doi-org.proxy.libraries.rutgers.edu/10.4007/annals.2015.182.1.1},
}

@article{LM15,
    AUTHOR = {Liu, Yi and Markovi\'c, Vladimir},
     TITLE = {Homology of curves and surfaces in closed hyperbolic 3-manifolds},
   JOURNAL = {Duke Math. J.},
  FJOURNAL = {Duke Mathematical Journal},
    VOLUME = {164},
      YEAR = {2015},
    NUMBER = {14},
     PAGES = {2723--2808},
      ISSN = {0012-7094,1547-7398},
   MRCLASS = {57M05 (20H10)},
  MRNUMBER = {3417184},
MRREVIEWER = {Mahan\ Mj},
       DOI = {10.1215/00127094-3167744},
       URL = {https://doi-org.proxy.libraries.rutgers.edu/10.1215/00127094-3167744},
}

@article {Sun15,
    AUTHOR = {Sun, Hongbin},
     TITLE = {Virtual homological torsion of closed hyperbolic 3-manifolds},
   JOURNAL = {J. Differential Geom.},
  FJOURNAL = {Journal of Differential Geometry},
    VOLUME = {100},
      YEAR = {2015},
    NUMBER = {3},
     PAGES = {547--583},
      ISSN = {0022-040X,1945-743X},
   MRCLASS = {57M27 (57M20)},
  MRNUMBER = {3352799},
MRREVIEWER = {Maki\ Nagura},
       URL = {http://projecteuclid.org/euclid.jdg/1432842365},
}

@article{Sun22,
    AUTHOR = {Sun, Hongbin},
     TITLE = {The panted cobordism groups of cusped hyperbolic 3-manifolds},
   JOURNAL = {J. Topol.},
  FJOURNAL = {Journal of Topology},
    VOLUME = {15},
      YEAR = {2022},
    NUMBER = {3},
     PAGES = {1580--1634},
      ISSN = {1753-8416,1753-8424},
   MRCLASS = {57K32 (57M05)},
  MRNUMBER = {4461854},
MRREVIEWER = {Ken'ichi\ Yoshida},
       DOI = {10.1112/topo.12255},
       URL = {https://doi-org.proxy.libraries.rutgers.edu/10.1112/topo.12255},
}

@article {Sun23,
    AUTHOR = {Sun, Hongbin},
     TITLE = {Virtual domination of 3-manifolds {II}},
   JOURNAL = {J. Lond. Math. Soc. (2)},
  FJOURNAL = {Journal of the London Mathematical Society. Second Series},
    VOLUME = {108},
      YEAR = {2023},
    NUMBER = {3},
     PAGES = {869--915},
      ISSN = {0024-6107,1469-7750},
   MRCLASS = {57K32 (30F40 57M10)},
  MRNUMBER = {4639942},
MRREVIEWER = {Colin\ C.\ Adams},
       DOI = {10.1112/jlms.12770},
       URL = {https://doi.org/10.1112/jlms.12770},
}

@misc{KW19,
  title={\emph{Counting connections in a locally symmetric space}},
  author={Kahn, Jeremy and Wright, Alex},
  howpublished = {\url{https://www.math.brown.edu/jk17/Connections.pdf}},
  year={2019}
}

@article {Rao25,
    AUTHOR = {Rao, Zhenghao},
     TITLE = {Subgroups of genus-2 quasi-{F}uchsian groups and cocompact
              {K}leinian groups},
   JOURNAL = {Geom. Topol.},
  FJOURNAL = {Geometry \& Topology},
    VOLUME = {29},
      YEAR = {2025},
    NUMBER = {1},
     PAGES = {495--548},
      ISSN = {1465-3060,1364-0380},
   MRCLASS = {57K32 (22E40 30C65 30F35 30F40 57K20)},
  MRNUMBER = {4877264},
MRREVIEWER = {Ken-ichi\ Ohshika},
       DOI = {10.2140/gt.2025.29.495},
       URL = {https://doi-org.proxy.libraries.rutgers.edu/10.2140/gt.2025.29.495},
}

@book {BP92,
    AUTHOR = {Benedetti, Riccardo and Petronio, Carlo},
     TITLE = {Lectures on hyperbolic geometry},
    SERIES = {Universitext},
 PUBLISHER = {Springer-Verlag, Berlin},
      YEAR = {1992},
     PAGES = {xiv+330},
      ISBN = {3-540-55534-X},
   MRCLASS = {57M50 (30F40 30F60 51M10 57N10)},
  MRNUMBER = {1219310},
MRREVIEWER = {Colin\ C.\ Adams},
       DOI = {10.1007/978-3-642-58158-8},
       URL = {https://doi-org.proxy.libraries.rutgers.edu/10.1007/978-3-642-58158-8},
}

@article {Agol,
    AUTHOR = {Agol, Ian},
     TITLE = {The virtual {H}aken conjecture},
      NOTE = {With an appendix by Agol, Daniel Groves, and Jason Manning},
   JOURNAL = {Doc. Math.},
  FJOURNAL = {Documenta Mathematica},
    VOLUME = {18},
      YEAR = {2013},
     PAGES = {1045--1087},
      ISSN = {1431-0635,1431-0643},
   MRCLASS = {20F67 (57Mxx)},
  MRNUMBER = {3104553},
MRREVIEWER = {Thomas\ Koberda},
       URL = {https://elibm.org/article/10000267},
}

@article{Ham15,
AUTHOR = {Hamenst\"adt, Ursula},
     TITLE = {Incompressible surfaces in rank one locally symmetric spaces},
   JOURNAL = {Geom. Funct. Anal.},
  FJOURNAL = {Geometric and Functional Analysis},
    VOLUME = {25},
      YEAR = {2015},
    NUMBER = {3},
     PAGES = {815--859},
      ISSN = {1016-443X,1420-8970},
   MRCLASS = {53C35 (22E40 37A25 37D40 53B21)},
  MRNUMBER = {3361773},
MRREVIEWER = {Bruno\ Duchesne},
       DOI = {10.1007/s00039-015-0330-y},
       URL = {https://doi-org.proxy.libraries.rutgers.edu/10.1007/s00039-015-0330-y},
}

@article{KR26,
    AUTHOR = {Kahn, Jeremy and Rao, Zhenghao},
     TITLE = {Surface subgroups for cocompact lattices of isometries of
              {$\mathbb H^{2n}$}},
   JOURNAL = {J. Topol.},
  FJOURNAL = {Journal of Topology},
    VOLUME = {19},
      YEAR = {2026},
    NUMBER = {3},
     PAGES = {Paper No. e70089, 45},
      ISSN = {1753-8416,1753-8424},
   MRCLASS = {22E40 (30F40)},
  MRNUMBER = {5103600},
       DOI = {10.1112/topo.70089},
       URL = {https://doi-org.proxy.libraries.rutgers.edu/10.1112/topo.70089},
}

@book {Fen89,
    AUTHOR = {Fenchel, Werner},
     TITLE = {Elementary geometry in hyperbolic space},
    SERIES = {De Gruyter Studies in Mathematics},
    VOLUME = {11},
      NOTE = {With an editorial by Heinz Bauer},
 PUBLISHER = {Walter de Gruyter \& Co., Berlin},
      YEAR = {1989},
     PAGES = {xii+225},
      ISBN = {3-11-011734-7},
   MRCLASS = {51M10},
  MRNUMBER = {1004006},
MRREVIEWER = {A.\ F.\ Beardon},
       DOI = {10.1515/9783110849455},
       URL = {https://doi-org.proxy.libraries.rutgers.edu/10.1515/9783110849455},
}

@article {Sun25,
    AUTHOR = {Sun, Hongbin},
     TITLE = {Virtual domination of 3-manifolds, {III}},
   JOURNAL = {Algebr. Geom. Topol.},
  FJOURNAL = {Algebraic \& Geometric Topology},
    VOLUME = {25},
      YEAR = {2025},
    NUMBER = {3},
     PAGES = {1599--1666},
      ISSN = {1472-2747,1472-2739},
   MRCLASS = {57K32 (30F40 57M10)},
  MRNUMBER = {4930573},
MRREVIEWER = {Colin\ C.\ Adams},
       DOI = {10.2140/agt.2025.25.1599},
       URL = {https://doi.org/10.2140/agt.2025.25.1599},
}

@article {Masters,
    AUTHOR = {Masters, Joseph D.},
     TITLE = {Kleinian groups with ubiquitous surface subgroups},
   JOURNAL = {Groups Geom. Dyn.},
  FJOURNAL = {Groups, Geometry, and Dynamics},
    VOLUME = {2},
      YEAR = {2008},
    NUMBER = {2},
     PAGES = {263--269},
      ISSN = {1661-7207,1661-7215},
   MRCLASS = {57M50 (20H10)},
  MRNUMBER = {2393182},
MRREVIEWER = {John\ G.\ Ratcliffe},
       DOI = {10.4171/GGD/39},
       URL = {https://doi-org.proxy.libraries.rutgers.edu/10.4171/GGD/39},
}

@article {KRS,
    AUTHOR = {Kolpakov, Alexander and Reid, Alan W. and Slavich, Leone},
     TITLE = {Embedding arithmetic hyperbolic manifolds},
   JOURNAL = {Math. Res. Lett.},
  FJOURNAL = {Mathematical Research Letters},
    VOLUME = {25},
      YEAR = {2018},
    NUMBER = {4},
     PAGES = {1305--1328},
      ISSN = {1073-2780,1945-001X},
   MRCLASS = {57N16},
  MRNUMBER = {3882165},
MRREVIEWER = {Michelle\ Chu},
       DOI = {10.4310/MRL.2018.v25.n4.a12},
       URL = {https://doi-org.proxy.libraries.rutgers.edu/10.4310/MRL.2018.v25.n4.a12},
}

@article {McReynolds,
    AUTHOR = {McReynolds, D. B.},
     TITLE = {Controlling manifold covers of orbifolds},
   JOURNAL = {Math. Res. Lett.},
  FJOURNAL = {Mathematical Research Letters},
    VOLUME = {16},
      YEAR = {2009},
    NUMBER = {4},
     PAGES = {651--662},
      ISSN = {1073-2780},
   MRCLASS = {57N16 (22F50 57M50)},
  MRNUMBER = {2525031},
MRREVIEWER = {Darren\ D.\ Long},
       DOI = {10.4310/MRL.2009.v16.n4.a8},
       URL = {https://doi-org.proxy.libraries.rutgers.edu/10.4310/MRL.2009.v16.n4.a8},
}

\end{document}